\documentclass[leqno,11pt,pdftex]{amsart}

\usepackage[top=25truemm,bottom=25truemm,left=25truemm,right=25truemm]{geometry}

\usepackage{amsmath,amsthm,amssymb,amsfonts,amscd,mathtools,color}
\usepackage[all]{xy} 
\usepackage{graphicx} 
\numberwithin{equation}{section}

\theoremstyle{plain} 
	\newtheorem{thm}{Theorem}[section]
	\newtheorem*{thm*}{Theorem}
	\newtheorem{cor}[thm]{Corollary}
	\newtheorem{lem}[thm]{Lemma}
	
	\newtheorem{prop}[thm]{Proposition}
	
	\newtheorem*{conj*}{Conjecture}
	
\theoremstyle{definition}
	\newtheorem{defn}[thm]{Definition}

\theoremstyle{remark}
	\newtheorem{rem}[thm]{Remark}
	
	\newtheorem*{pf}{Proof}
\def\CC{{\mathbb C}}

\def\MM{{\mathbb M}}

\def\PP{{\mathbb P}}
\def\QQ{{\mathbb Q}}
\def\RR{{\mathbb R}}
\def\SS{{\mathbb S}}
\def\TT{{\mathbb T}}

\def\ZZ{{\mathbb Z}}
\def\A{{\mathcal A}}

\def\C{{\mathcal C}}
\def\D{{\mathcal D}}
\def\E{{\mathcal E}}
\def\F{{\mathcal F}}

\def\H{{\mathcal H}}

\def\L{{\mathcal L}}
\def\M{{\mathcal M}}

\def\O{{\mathcal O}}

\def\R{{\mathcal R}}
\def\S{{\mathcal S}}
\def\T{{\mathcal T}}
\def\U{{\mathcal U}}

\def\W{{\mathcal W}}

\def\p{{\partial}}

\def\Aut{{\rm Aut}}
\def\ST{{\rm ST}}
\def\coh{{\rm coh}}

\def\Exc{{\rm Exc}}
\def\Ext{{\rm Ext}}
\def\FEC{{\rm FEC}}
\def\Fuk{{\rm Fuk}}
\def\Hom{{\rm Hom}}
\def\Ind{{\rm Ind}}

\def\Sph{{\rm Sph}}
\def\Stab{{\rm Stab}}
\def\mod{{\rm mod}}

\def\h{{\mathfrak h}}

\def\vd{\vec{\Delta}}

\def\ns{{\nabla}\hspace{-1.4mm}\raisebox{0.3mm}{\text{\footnotesize{\bf /}}}}

\begin{document}
\title[The numbers of full exceptional collections for extended Dynkin quivers]{The number of full exceptional collections modulo spherical twists
for extended Dynkin quivers}
\author{Takumi Otani}
\address{Yau Mathematical Sciences Center, Tsinghua University, 
Haidian District, Beijing, 100084, China}
\email{takumi@tsinghua.edu.cn}
\author{Yuuki Shiraishi}
\address{School of Economics and Management, University of Hyogo, 
Nishiku Kobe Hyogo, 651-2197, Japan}
\email{s912y025@guh.u-hyogo.ac.jp}
\author{Atsushi Takahashi}
\address{Department of Mathematics, Graduate School of Science, The University of Osaka, 
Toyonaka Osaka, 560-0043, Japan}
\email{takahashi@math.sci.osaka-u.ac.jp}

\thanks{
© 2026. 
This manuscript version is made available under the CC-BY-NC-ND 4.0 license 
https://creativecommons.org/licenses/by-nc-nd/4.0/ \\
The Version of Record is available at https://doi.org/10.1016/j.aim.2026.111201.
}
\maketitle
\begin{abstract}
This paper calculates the number of full exceptional collections modulo an action of a free abelian group of rank one
for an abelian category of coherent sheaves on an orbifold projective line with a positive orbifold Euler characteristic,
which is equivalent to the one of finite dimensional modules over an extended Dynkin quiver of ADE type by taking their derived categories.
This is done by a recursive formula naturally generalizing the one for the Dynkin case by Deligne whose categorical interpretation is due to Obaid--Nauman--Shammakh--Fakieh--Ringel.

Moreover, the number coincides with the degree of the Lyashko--Looijenga map of the Frobenius manifold for the orbifold projective line, which hints a consistency in some problems in Bridgeland's stability conditions and mirror symmetry.
\end{abstract}
\section{Introduction}

Inspired by the correspondence between maximal chains in the poset of noncrossing partitions for a Dynkin quiver $\vd$
and complete exceptional sequences in the derived category $\D^b(\CC\vd)$ of finitely generated modules over 
the path algebra $\CC\vd$, Obaid--Nauman-- Shammakh--Fakieh--Ringel \cite{ONSFR} found the following recursive formula
for $e(\D^b(\CC\vd)) \coloneqq |\FEC(\D^b(\CC\vd)) /  \ZZ^\mu |$
where $\FEC(\D^b(\CC\vd))$ is the set of isomorphism classes of full exceptional collections in $\D^b(\CC\vd)$,
$\mu$ is the number of vertices of $\vd$ and $\ZZ^\mu$ is the group whose $i$-th generator acts by the translation functor $[1]$ on the $i$-th object of a full exceptional collection:
\begin{equation*}
e(\D^b(\CC\vd)) = \dfrac{h}{2} \sum_{v \in \Delta_0} e(\D^b(\CC\vd^{(v)})),
\end{equation*}
where $\vd^{(v)}$ is the full subquiver of $\vd$ given by removing the vertex $v$ and arrows connecting with $v$ and 
$h \in \ZZ_{\ge 0}$ is the Coxeter number. As is mentioned in \cite{ONSFR}, the above recursive formula 
is a categorification of Deligne's recursive formula \cite{De} which originally answered to Looijenga's conjecture for simple singularities concerning 
the number of distinguished bases of vanishing cycles modulo signs and the degree of a map $LL$, called the Lyashko--Looijenga 
map, describing the topology of the bifurcation set \cite{L}. 
As a summary of the above story, the following equation holds: 
\begin{equation*}
e(\D^b(\CC\vd)) = \frac{\mu!}{d_1\cdots d_\mu} h^\mu = \deg LL,
\end{equation*}
where $2\le d_1\le\dots\le d_\mu=h$ are degrees of algebraic independent invariants of the Weyl group for the Dynkin diagram $\Delta$.
Here $LL$ was originally defined by using the characteristic polynomial associated to 
the discriminant for the singularity, or equivalently, the square of Jacobian for the  Weyl group invariant theory.
More generally, $LL$ can be defined for massive $F$-manifolds (see Section \ref{sec : generalized LL map} and \cite{DZ, He}).
In \cite{HR}, Hertling--Roucairol extended the above correspondence between the number of distinguished bases 
of vanishing cycles modulo signs and the degree of the Lyashko--Looijenga map to simple elliptic singularities with Legendre normal forms.
Therefore, it is very important to understand what happens when the setting falls 
between the two classes, corresponding to simple singularities and the simple elliptic singularities.

The purpose of this paper is to generalize the above recursive formula and equality for the derived category $\D^b(\PP_A^1)$ of
coherent sheaves over an orbifold projective line $\PP_A^1$ with three orbifold points $(0,1,\infty)$ whose orders are given by 
$A \coloneqq (a_{1}, a_{2}, a_{3})$ satisfying $\displaystyle \chi_{A} \coloneqq \frac{1}{a_{1}}+\frac{1}{a_{2}}+\frac{1}{a_{3}}-1> 0$. 
Note that the condition $\chi_{A}>0$ implies that $\D^b(\PP_A^1)$ is equivalent to the derived category $\D^{b}(\CC Q_A)$ 
of an extended Dynkin quiver $Q_A$ due to Geigle--Lenzing \cite{GL1} (see Proposition \ref{prop : derived equivalence}). 
The one of most crucial points of this paper is to define $e(\D^b(\PP_A^1))$ from the perspective of mirror symmetry
and the space of stability conditions. It is conjectured and is proved for $A=(1,p,q)$ in \cite{HKK} that 
\begin{equation*}
M_{\PP_A^1}=\CC^{\mu_{A}-1}\times \CC^{*}\cong \Stab(\D^b(\PP_A^1))/\ZZ=\CC^{\mu_{A}}/\ZZ,
\end{equation*}
where $M_{\PP_A^1}$ is the Frobenius manifold of rank $\mu_{A} \coloneqq a_{1}+a_{2}+a_{3}-1$ from the orbifold Gromov--Witten theory for 
$\PP_A^1$, $\Stab(\D^b(\PP_A^1))$ is the space of stability conditions for $\D^b(\PP_A^1)$ and $\ZZ$ is a subgroup of 
the autoequivalence group $\Aut(\D^b(\PP_A^1))$ of $\D^b(\PP_A^1)$. 
Especially, we conjecture that this $\ZZ$ is the subgroup $\ST(\D^b(\PP_A^1))$ generated by the spherical twist
$-\otimes \O(\vec{c})$ (see Proposition~\ref{prop : Braid group of an orbifold projective line}). 
Due to the polynomiality of the Frobenius potential for $M_{\PP_A^1}$, the degree of the Lyashko--Looijenga map  makes sense.
Thus we reach to our definition of $e(\D^b(\PP_A^1))$ as 
\begin{equation*}
e(\D^b(\PP_A^1)) \coloneqq |\FEC(\D^b(\PP_A^1)) / \langle \ST(\D^b(\PP_A^1)), \ZZ^{\mu_A} \rangle|.
\end{equation*}

Under the above conjecture and our definition,
the first main result is a Deligne type recursive formula of the number $e(\D^b(\PP_A^1))$.
For each vertex $v \in Q_A$, let $Q_A^{(v)}$ denote the full subquiver of the extended Dynkin quiver $Q_A$ whose vertices are given by $(Q_A)_0 \setminus \{ v \}$. Note that the quiver $Q_A^{(v)}$ is a union of Dynkin quivers.
Let $\vec{A}_{a_i - j - 1}$ denote the Dynkin quiver of type $A_{a_i - j - 1}$ with all arrows oriented to the left. 
\begin{thm}[Theorem \ref{thm : recursion formula}]\label{thm : intro recursion formula}
We have
\begin{equation*}
e(\D^b(\PP_A^1)) = \dfrac{1}{\chi_A} \sum_{v \in (Q_A)_0} e(\D^b(\CC Q_A^{(v)})) + \sum_{i = 1}^3 a_i \sum_{j = 1}^{a_i - 1}\binom{\mu_A-1}{a_i - j-1}\cdot e(\D^b(\PP^1_{A_{(i, j)}}))\cdot e( \D^b(\CC \vec{A}_{a_i - j - 1})),
\end{equation*}
where $A_{(i, j)} = (a_1', a_2', a_3')$ is defined by $a'_i = j$ and $a'_k = a_k$ for $k \ne i$.
\end{thm}
For $M_{\PP_A^1}$, the degree of  the Lyashko--Looijenga map is calculated in \cite[Section 3]{DZ} as
\begin{equation*}
\deg LL = \frac{\mu_A!}{a_1!a_2!a_3!\chi_A}a_1^{a_1}a_2^{a_2}a_3^{a_3}.
\end{equation*}
The second main result is the equality of $e(\D^b(\PP_A^1))$ and the degree of the Lyashko--Looijenga map:
\begin{thm}[Theorem \ref{thm : main}, Corollary \ref{cor: equality to deg LL}]\label{thm : intro main}
Let $\PP_A^1$ be an orbifold projective line with a positive orbifold Euler characteristic.
The number $e(\D^b(\PP_A^1))$ is equal to the degree of the Lyashko--Looijenga map of the Frobenius manifold $M_{\PP_A^1}$: 
\begin{equation*}
e(\D^b(\PP_A^1)) = \deg LL .
\end{equation*}
\end{thm}
The above equality when $A=(1,p,q), (2,2,r)$, corresponding to the extended Dynkin diagram of type A and D, can be shown by recursive formula due to Hurwitz~\cite{Hu} (see also \cite{S-V}). 
Indeed, this formula and the observation of $\ST(\D^b(\PP_A^1))$ as a mapping class group match with the proof in \cite{HKK} 
using the marked bordered annulus, which inspire future directions to generalize our work to 
the partially wrapped Fukaya category $\W(\SS, \MM; \ell)$ and 
a Frobenius structure associated to a tuple $(\SS, \MM; \ell)$.
Here $\SS$ is a surface with boundaries, 
$\MM$ is the disjoint union consisting of $k_{i}$ marked points on each boundary $\p_{i}$ of $\SS$ for $i=1,\dots b$ and   
$\ell$ is a particular choice of a line field.
More precisely, when $\D=\W(\SS, \MM; \ell)$,
the mapping class group ${\rm MCG}(\SS)$ acts on $\FEC(\D)$.
We expect that $\ST(\D)$ should be isomorphic to ${\rm MCG}(\SS)$.
The results in \cite{CHS, ELSV, Klu} yield the following formula:
\begin{equation*}
e(\D) \coloneqq |\FEC(\D) / \langle {\rm MCG}(\SS), \ZZ^{\mu} \rangle|=k_{1}\cdots k_{b}\cdot
|{\rm Aut}(k_{1},\dots, k_{b})|\cdot h_{g;k_{1},\dots, k_{b}}=k_{1}\cdots k_{b} \cdot \deg \L\L,
\end{equation*}
where we use the symbols in \cite{ELSV}, i.e.,
$|{\rm Aut}(k_{1},\dots, k_{b})|$ is the number of permutations preserving the sequence $(k_{1},\dots, k_{b})$, 
$h_{g;k_{1},\dots, k_{b}}$ and $\deg \L\L$ are the Hurwitz number and the degree of the Lyashko--Looijenga map 
for the Hurwitz space. 
As a result, the number $e(\D)$ is equal to the degree of a Lyashko--Looijenga map for a space of meromorphic functions with marked branches of roots at poles, a $\ZZ/k_{1}\ZZ\times\dots \times \ZZ/k_{b}\ZZ$--cover of the Hurwitz space, on which 
Dubrovin considered a (twisted) Frobenius structure. 
Moreover, the Deligne type recursive relation for $e(\D)$ coincides with the cut-and-join equation by \cite{GJ}.
As another explanation for the importance of $e(\D)$,
we also note that a marked bordered annulus with a choice of line fields plays
important roles not only as the mirror partner to $\PP^{1}_{A}$ with $A=(1,p,q)$ by the line field
corresponding to the primitive form in \cite{MT}, but also for constructing,  
via a marked bordered annulus with an involution and another line field corresponding to the different one in \cite{T2},
a root system of type $D$ with a quasi-Coxeter element geometrically and a Frobenius structure associated to it. 
Then, the above definition of $e(\D)$ is also compatible with this construction  
from the viewpoint of the degree of the Lyashko--Looijenga map and counting the real roots.
The details of the above results will be written elsewhere \cite{OIST}.

Recently, similar results are given for simple elliptic singularities with Legendre normal forms \cite{TZ}.
In particular, \cite{TZ} gives recursive expressions of the Lyashko--Looijenga maps in terms of  $A$ based on
a cut-and-join equation in the categorical setting, while the ones by \cite{HR} is based on weights and 
can be considered as a variation of ELSV formulas.

Therefore, Seidel--Thomas groups have the potential to understand Frobenius manifolds globally and
to make a precise connection to the spaces of stability conditions via the set of full exceptional collections.

\bigskip
\noindent
{\bf Acknowledgements.}
Through joint works with him, the authors are grateful to Akishi Ikeda for valuable discussions
on spaces of stability conditions associated to marked bordered surfaces, Hurwitz--Frobenius structures
and generalized root systems.    
T.O and A.T are grateful to Claus Hertling for valuable discussions about \cite{HR}, which is one motivation for this paper.
A.T is supported by JSPS KAKENHI Grant Number JP21H04994.
This research and the paper were completed during his affiliation of T. O was Osaka University.
T. O was partially supported by JSPS KAKENHI Grant Number JP21H04994.
Y.S is supported by JSPS KAKENHI Grant Numbers 19K14531 and 23K03111.
\bigskip
\noindent
{\bf Notation.}
Throughout this paper, for a finite-dimensional $\CC$-algebra $A$, the bounded derived category of finitely generated right $A$-modules 
is denoted by $\D^b(A) \coloneqq \D^b \mod (A)$.
We also denote by $\D^b(X)$ the derived category $\D^b \coh(X)$ of coherent sheaves on an orbifold $X$.
For a triangulated category $\D$, the group of autoequivalences of $\D$ is denoted by $\Aut(\D)$.

\section{Preliminaries}\label{sec : preliminaries}
\subsection{Full exceptional collections}
First, we recall the notion of a full exceptional collection and basic properties.
\begin{defn}
Let $\D$ be a $\CC$-linear triangulated category.
An object $E \in \D$ is called {\em exceptional} if $\Hom_\D(E,E) \cong \CC$ and $\Hom_\D(E,E[p]) \cong 0$ when $p \ne 0$.
\end{defn}
An object $E \in \D$ is said to be {\em indecomposable} if $E$ is nonzero and $E$ has no direct sum decomposition $E \cong A \oplus B$, where $A$ and $B$ are nonzero objects in $\D$.
Note that an exceptional object is indecomposable.
Denote by $\Exc(\D)$ (resp. $\Ind(\D)$) the set of isomorphism classes of exceptional (resp. indecomposable) objects in $\D$.
For each exceptional object $E \in \D$, we define a full subcategory $E^\perp$ by 
\begin{equation*}
E^\perp \coloneqq \{ X \in \D \mid \Hom_\D(E, X[p]) = 0 ~ \text{for all} ~ p \in \ZZ \},
\end{equation*}
which gives a semi-orthogonal decomposition $\D = \langle E^\perp, E \rangle$.

\begin{defn}
Let $\D$ be a $\CC$-linear triangulated category.
\begin{enumerate}
\item An ordered set $\E = (E_1, \dots, E_\mu)$ consisting of exceptional objects $E_1, \dots, E_\mu$ is called {\em exceptional collection} if $\Hom^p_\D(E_i, E_j) \cong 0$ for all $p \in \ZZ$ and $i > j$.
\item An exceptional collection $\E$ is called {\em full} if the smallest full triangulated subcategory of $\D$ containing all elements in $\E$ is equivalent to $\D$ as a triangulated category.
\item Two full exceptional collections $\E=(E_1, \dots, E_\mu)$, $\F=(F_1, \dots, F_\mu)$ in $\D$ are said to be {\em isomorphic} if $E_i \cong F_i$ for all $i = 1, \dots, \mu$.
\end{enumerate}
\end{defn}

\begin{defn}
Let $(E, F)$ be an exceptional collection.
Define two objects $\R_F E$ and $\L_E F$ by the following exact triangles respectively:
\begin{equation*}
\R_F E \longrightarrow E \overset{{\rm ev}^*}{\longrightarrow} \Hom^\bullet_\D(E, F)^* \otimes F,
\end{equation*}
\begin{equation*}
\Hom^\bullet_\D(E, F) \otimes E \overset{\rm ev}{\longrightarrow} F \longrightarrow \L_E F,
\end{equation*}
where $(-)^*$ denotes the duality $\Hom_\CC(-, \CC)$.
The object $\R_F E$ (resp. $\L_E F$) is called the {\em right mutation} of $E$ through $F$ (resp. {\em left mutation} of $F$ through $E$).
Then, $(F, \R_F E)$ and $(\L_E F, E)$ form new exceptional collections.
\end{defn}
\begin{rem}
Our definition of mutations differs from the usual one (cf. \cite[Section 1]{BP}).
In our notation, the usual left and right mutations are given by $\R_F E [1]$ and $\L_E F [-1]$, respectively.
\end{rem}

The Artin's {\it braid group} ${\rm Br}_{\mu}$ on $\mu$-stands is a group presented by the following generators and relations: 
\begin{description}
\item[{\bf Generators}] $\{b_i~|~i=1,\dots, \mu-1\}$
\item[{\bf Relations}] $b_{i}b_{j}=b_{j}b_{i}$ for $|i-j|\ge 2$, $b_{i}b_{i+1}b_{i}=b_{i+1}b_{i}b_{i+1}$ for $i=1,\dots, \mu-2$.
\end{description}
Consider the group ${\rm Br}_{\mu} \ltimes \ZZ^{\mu}$, the semi-direct product of the braid group ${\rm Br}_{\mu}$ and the 
abelian group $\ZZ^{\mu}$, defined by the group homomorphism ${\rm Br}_{\mu} \longrightarrow {\mathfrak S}_{\mu} \longrightarrow {\rm Aut}_{\ZZ} \ZZ^{\mu}$, where the first homomorphism is $b_{i}\mapsto (i, i+1)$ and the second one is induced by the natural actions of the symmetric group ${\mathfrak S}_{\mu}$ on $\ZZ^{\mu}$. 
\begin{prop}[{cf.~\cite[Proposition 2.1]{BP}}]\label{prop : braid}
Let $\D$ be a $\CC$-linear triangulated category. Define $\FEC(\D)$ as the set of isomorphism classes of full exceptional collections in $\D$.
The group ${\rm Br}_{\mu} \ltimes \ZZ^{\mu}$ acts on $\FEC(\D)$ by mutations and transformations:
\begin{gather*}
b_{i} \cdot (E_{1}, \dots, E_{\mu}) \coloneqq (E_{1}, \dots, E_{i-1}, E_{i+1}, \R_{E_{i+1}} E_i, E_{i+2}, \dots, E_{\mu}), \\
b^{-1}_{i} \cdot (E_{1}, \dots, E_{\mu}) \coloneqq (E_{1}, \dots, E_{i-1}, \L_{E_{i}} E_{i+1}, E_{i}, E_{i+2},\dots, E_{\mu}), \\
e_{i} \cdot (E_{1}, \dots, E_{\mu}) \coloneqq (E_{1}, \dots, E_{i-1}, E_{i} [1], E_{i+1}, \dots, E_{\mu}),
\end{gather*}
where $e_{i}$ is the $i$-th generator of $\ZZ^{\mu}$.
\qed
\end{prop}

In this paper, we shall consider the derived category $\D^b(\A)$ of a {\em hereditary} abelian category $\A$, 
namely, an abelian category $\A$ satisfying $\Ext^p_\A(X, Y) = 0$ for any $p \ge 2$ and objects $X, Y \in \A$.
Indecomposable objects and exceptional objects in such derived category have the property that 
\begin{equation*}
\Ind (\D^b(\A)) \cong \bigsqcup_{p \in \ZZ} \Ind (\A) [p], \quad 
\Exc (\D^b(\A)) \cong \bigsqcup_{p \in \ZZ} \Exc (\A) [p],
\end{equation*}
where $\Ind (\A) [p]$ (resp. $\Exc (\A) [p]$) denotes the set of isomorphism classes of indecomposable (resp. exceptional) objects $E \in \A$ translated by $p$.

For a $\CC$-linear triangulated category $\D$, denote by $\Aut(\D)$ the group of auto-equivalences.
Finally, we introduce an action of $\Aut(\D)$ on $\FEC(\D)$.
\begin{defn}
Let $\D$ be a $\CC$-linear triangulated category. 
The group of auto-equivalences of $\D$ acts on $\FEC(\D)$ by 
\[
\Phi \cdot (E_1, \dots, E_\mu) \coloneqq (\Phi(E_1), \dots, \Phi(E_\mu)), \quad \Phi \in \Aut(\D).
\]
\end{defn}

\section{Orbifold projective lines}\label{sec : orbifold projective line}
We summarize the notations and results in \cite{GL1}.
Let $A = (a_1, a_2 ,a_3)$ be a triplet of positive integers. 
Let $a$ denote the least common multiple of $a_1, a_2$ and $a_3$.

\begin{defn}[{\cite[Section 1.1]{GL1}}]
Let $A = (a_1, a_2, a_3)$.
\begin{enumerate}
\item Denote by $L_A$ an abelian group generated by $3$-letters $\vec{x_i}$, $i = 1, 2, 3$ defined as the quotient 
\begin{equation*}
L_A \coloneqq \left. \bigoplus_{i = 1}^3 \ZZ \vec{x}_i \middle/ \big\langle a_i \vec{x}_i - a_j \vec{x}_j \mid i, j = 1, 2, 3 \big\rangle \right. .
\end{equation*}
We define the {\em canonical element} of $L_A$ by $\vec{c} \coloneqq a_1 \vec{x}_1 = a_2 \vec{x}_2 = a_3 \vec{x}_3$.
\item Define an $L_A$-graded $\CC$-algebra $S_{A}$ by 
\begin{equation*} 
S_{A} \coloneqq \CC[X_1, X_2, X_3] / (X_3^{a_3} - X_2^{a_2} + X_1^{a_1}) ,
\end{equation*}
where $\deg (X_i) = \vec{x}_i$ for $i = 1, 2, 3$.
\end{enumerate}
\end{defn}

\begin{defn}[{cf.~\cite[Section 1.1]{GL1}}]
Let $A = (a_1, a_2, a_3)$.
Define a stack $\PP^1_{A}$ by
\begin{equation*}
\PP_A^1 \coloneqq \left[ \left( {\rm Spec} (S_A) \backslash \{ 0 \} \right) / {\rm Spec} ({\CC L_A}) \right].
\end{equation*}
The stack $\PP_A^1$ is called the {\em orbifold projective line}. 
\end{defn}

Denote by $\mod^{L_A}(S_{A})$ the abelian category of finitely generated $L_A$-graded $S_{A}$-modules and denote by $\mod^{L_A}_0 (S_{A})$ the full subcategory of ${\rm gr}^{L_A}(S_{A})$ whose objects are $L_A$-graded finite length $S_{A}$-modules.
It is known by \cite[Section 1.8]{GL1} that the abelian category $\coh(\PP^{1}_{A})$ of coherent sheaves on $\PP^{1}_{A}$ is given by 
\[
\coh(\PP^{1}_{A}) = \mod^{L_A} (S_{A}) / \mod^{L_A}_0 (S_{A}),
\]

\begin{prop}[{\cite[Section 2.2]{GL1}}]\label{prop : hereditary}
The abelian category $\coh(\PP_A^1)$ is hereditary.
\qed
\end{prop}

Define a sheaf $\O (\vec{l})$ for $\vec{l} \in L_{A}$ by
\begin{equation*}
\O (\vec{l}) \coloneqq [ S_{A} (\vec{l}) ] \in \coh(\PP_A^1)
\end{equation*}
where the degree $\vec{l'}$ homogeneous component $(S_{A}(\vec{l}))_{\vec{l'}}$ of $S_{A}(\vec{l})$
is defined by the degree $\vec{l} + \vec{l'}$ homogeneous component $(S_{A})_{\vec{l} + \vec{l'}}$ of $S_{A}$. 

\begin{defn}[{\cite[Section 2.5]{GL1}}]
Take $\lambda \in \PP^{1} \setminus \{ 0, 1, \infty \}$. 
Define $S_\lambda$ and $S_{i, j}$ for $i = 1, 2, 3$ and $j = 0, \dots, a_{i} - 1$ by the following exact sequences:
\begin{subequations}
\begin{equation}\label{eq : simple at an ordinary point}
0\rightarrow \O \xrightarrow{X^{a_{2}}_{2} - \lambda X^{a_{1}}_{1}} \O (\vec{c}) \rightarrow S_\lambda \rightarrow 0,
\end{equation}
\begin{equation}
0\rightarrow \O (j \vec{x}_{i}) \xrightarrow{X_{i}} \O ((j + 1) \vec{x}_{i}) \rightarrow S_{i, j} \rightarrow 0.
\end{equation}
\end{subequations}
\end{defn}

A coherent sheaf on $\PP_A^1$ is said to be a {\em vector bundle}  if the sheaf is locally free.
We denote by ${\rm vect}(\PP_A^1)$ the full subcategory of $\coh(\PP_A^1)$ consisting of vector bundles on $\PP_A^1$.
Each coherent sheaf $\F \in \coh(\PP_A^1)$ splits into a direct sum $\F_+ \oplus \F_0$ of $\F_+ \in {\rm vect} (\PP_A^1)$ and $\F_0 \in \coh_0 (\PP_A^1)$, where $\coh_0(\PP_A^1)$ is the full subcategory consisting of coherent sheaves of finite length.
The abelian category $\coh_0 (\PP_A^1)$ decomposes into a coproduct $\coprod_{\lambda \in \PP^1} \U_\lambda$, where $\U_\lambda$ denotes the uniserial category of finite length sheaves concentrated at the point $\lambda$ (\cite[Proposition 2.4 and 2.5]{GL1}).
Here an abelian category is said to be uniserial if the sub-objects of any indecomposable object form a chain under inclusion.

\begin{prop}[{\cite[Section 2.2]{GL1}}]\label{prop : Serre functor}
Define an element $\vec{\omega} \in L_A$ by 
\begin{equation}
\vec{\omega} \coloneqq \vec{c} - \sum_{i =  1}^3 \vec{x}_i.
\end{equation}
Then, the autoequivalence $\S_{\D^b(\PP_A^1)} \coloneqq - \otimes \O (\vec{\omega}) [1] \in \Aut \, \D^b (\PP_A^1)$ is the Serre functor.
Namely, there is an isomorphism
\[
\Hom_{\D^b(\PP_A^1)}(X, Y) \cong \Hom_{\D^b(\PP_A^1)}(Y, \S_{\D^b(\PP_A^1)}(X))^*,
\]
which is functorial with respect to $X \in \D^b(\PP_A^1)$ and $Y \in \D^b(\PP_A^1)$.
\qed
\end{prop}
When the category $\D^b(\PP_A^1)$ is clear from the context, we often drop it from the notation and write $\S$.

Define a positive integer $\mu_A$ by 
\begin{equation*}
\mu_A \coloneqq 2 + \sum_{i = 1}^3 (a_i - 1).
\end{equation*}

\begin{prop}[{\cite[Section 4.1]{GL1}}]\label{prop : full strongly exceptional collection}
The ordered set of $\mu_A$ exceptional sheaves on $\PP_A^1$
\[
( \O, \O(\vec{x}_1), \cdots, \O((a_1 - 1)\vec{x}_1)), \O(\vec{x}_2), \cdots, \O((a_2 - 1)\vec{x}_2)), \O(\vec{x}_3), \cdots, \O((a_3 - 1)\vec{x}_3), \O(\vec{c}) )
\]
is a full strongly exceptional collection in $\D^b(\PP_{A}^1)$.
\qed
\end{prop}

Define a rational number $\chi_A$ by 
\begin{equation*}
\chi_A \coloneqq 2 + \sum_{i = 1}^{3} \left( \frac{1}{a_i} - 1 \right),
\end{equation*}
which is called the {\em orbifold Euler characteristic} of $\PP_A^1$.
An orbifold projective line $\PP_A^1$ satisfying $\chi_A > 0$ is often called of {\em domestic type}.
It is important to note that $\chi_A>0$ with $a_1 \le a_{2}\le a_{3}$
if and only if $A = (1, p, q), (2, 2, r), (2, 3, 3), (2, 3, 4)$ or $(2, 3, 5)$, where $p, q, r \in \ZZ_{\ge 1}$.

Throughout this paper, we only treat orbifold projective lines with a positive orbifold Euler characteristic.

\subsection{Derived equivalence}

\begin{defn}[{\cite[Definition 2.19]{STW}}]
Define a quiver with relation $\widetilde{\TT}_A = ((\widetilde{\TT}_A)_0, (\widetilde{\TT}_A)_1, I)$ as follows:
\begin{itemize}
\setlength{\leftskip}{-10pt}
\item The set of vertices is given by 
\[
(\widetilde{\TT}_A)_0 \coloneqq \{ (i, j) \mid i = 1, 2, 3, ~ j = 1, \dots, a_i - 1 \} \sqcup \{ \mathbf{1}, \mathbf{1}^* \}.
\]
\item Let $v, v' \in (\widetilde{\TT}_A)_0$ be verticies.
\begin{itemize}
\setlength{\leftskip}{-10pt}
\item If $(v, v') = (\mathbf{1}, (i, 1))$ or $((i, 1), \mathbf{1}^*)$, there is one arrow $\alpha_{v, v'} \in (\widetilde{\TT}_A)_1$ from $v$ to $v'$.
\item If $(v, v') = ((i, j), (i, j + 1))$ for some $i, j$, there is one arrow $\alpha_{v, v'} \in (\widetilde{\TT}_A)_1$ from $v$ to $v'$.
\item Otherwise, there are no arrows.
\end{itemize}\item The relation $I$ is given by 
\[
I \coloneqq \langle \alpha_{(1, 1), \mathbf{1}^*} \alpha_{\mathbf{1}, (1, 1)} + \alpha_{(2, 1), \mathbf{1}^*} \alpha_{\mathbf{1}, (2, 1)}, \alpha_{(2, 1), \mathbf{1}^*} \alpha_{\mathbf{1}, (2, 1)} + \alpha_{(3, 1), \mathbf{1}^*} \alpha_{\mathbf{1}, (3, 1)} \rangle.
\]
\end{itemize}
\end{defn}
Figure \ref{fig : octopus} shows the quiver with relation $\widetilde{\TT}_A$.
For simplicity, we denote by $\CC \widetilde{\TT}_A$ the bounded quiver algebra $\CC \big( (\widetilde{\TT}_A)_0, (\widetilde{\TT}_A)_1 \big) / I$. 

\begin{figure}[h]
\[
\widetilde{\TT}_{A} \colon 
\xymatrix{
& & & \overset{\mathbf{1}^*}{\circ} \ar@{==}[d] & & & \\
\underset{(1, a_1 - 1)}{\circ} & \cdots \ar[l] & \underset{(1, 1)}{\circ} \ar[l] \ar[ru] & \underset{\mathbf{1}}{\circ} \ar[l] \ar[r] \ar[ld] & \underset{(3, 1)}{\circ} \ar[r] \ar[lu] & \cdots \ar[r] & \underset{(3, a_3 - 1)}{\circ} \\
& & \underset{(2, 1)}{\circ} \ar[ruu] \ar[ld] & & & & \\
& \rotatebox{75}{$\ddots$} \ar[ld] & & & & & \\
\underset{(2, a_2 - 1)}{\circ} & & & & & & 
}
\]
\caption{The quiver with relation $\widetilde{\TT}_A$.}
\label{fig : octopus}
\end{figure}
For each vertex $(i, j) \in (\widetilde{\TT}_A)_0$, we denote by $E^{\widetilde{\TT}_{A}}_{(i, j)} \in \D^b(\CC \widetilde{\TT}_A)$ the exceptional object corresponding to 
the projective $\CC \widetilde{\TT}_A$-module with respect to the vertex $(i, j) \in (\widetilde{\TT}_{A})_0$.

Next, we consider extended Dynkin quivers.
An extended Dynkin quiver is one of the following acyclic quivers:
\begin{itemize}
\item $A_{p, q}^{(1)}$-quiver: 
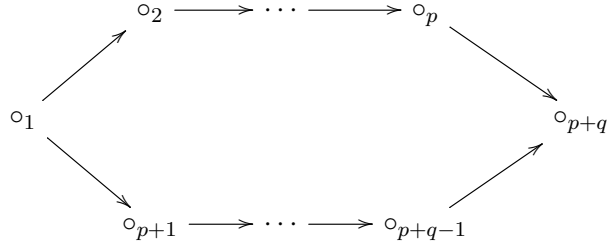
\begin{equation*}
\xymatrix{
& \circ_2 \ar[r] & \cdots \ar[r] & \circ_p \ar[rd] & \\
\circ_1 \ar[ru] \ar[rd] & & & & \circ_{p + q} \\
& \circ_{p + 1} \ar[r] & \cdots \ar[r] & \circ_{p + q - 1} \ar[ru] & 
}
\end{equation*}
\item $D_{r}^{(1)}$-quiver: 
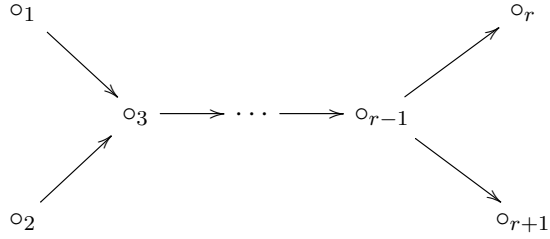
\begin{equation*}
\xymatrix{
\circ_1 \ar[rd] & & & & \circ_{r} \\
& \circ_3 \ar[r] & \cdots \ar[r] & \circ_{r - 1} \ar[ru] \ar[rd] & \\
\circ_2 \ar[ru] & & & & \circ_{r + 1}
}
\end{equation*}
\item $E_{6}^{(1)}$-quiver: 
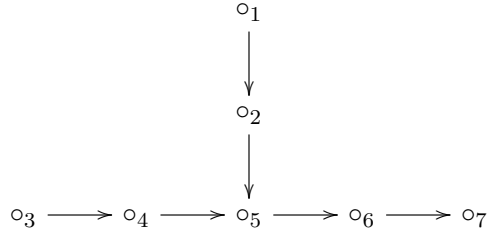
\begin{equation*}
\xymatrix{
& & \circ_1 \ar[d] & & \\
& & \circ_2 \ar[d] & & \\
\circ_3 \ar[r] & \circ_4 \ar[r] & \circ_5 \ar[r] & \circ_6 \ar[r] & \circ_7
}
\end{equation*}
\item $E_{7}^{(1)}$-quiver: 
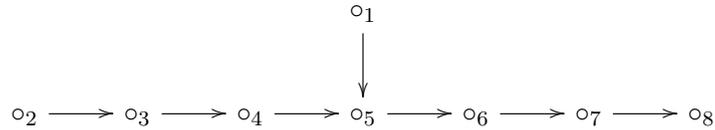
\begin{equation*}
\xymatrix{
& & & \circ_1 \ar[d] & & & \\
\circ_2 \ar[r] & \circ_3 \ar[r] & \circ_4 \ar[r] & \circ_5 \ar[r] & \circ_6 \ar[r] & \circ_7 \ar[r] & \circ_8
}
\end{equation*}
\item $E_{8}^{(1)}$-quiver: 
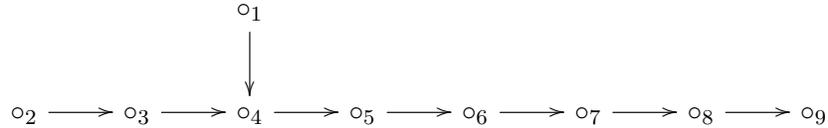
\begin{equation*}
\xymatrix{
& & \circ_1 \ar[d] & & & & & \\
\circ_2 \ar[r] & \circ_3 \ar[r] & \circ_4 \ar[r] & \circ_5 \ar[r] & \circ_6 \ar[r] & \circ_7 \ar[r] & \circ_8 \ar[r] & \circ_9
}
\end{equation*}
\end{itemize}
For a vertex $v \in Q_0$ of an extended Dynkin quiver $Q$, let $E^{Q}_v$ denote the exceptional object corresponding to 
the projective $\CC Q$-module with respect to the vertex $v \in Q_0$. 

For each vertex $v \in Q_0$, we define a quiver $Q^{(v)}$ as the full subquiver of $Q$ whose vertices are given by $Q_0 \setminus \{ v \}$.
Note that the quiver $Q^{(v)}$ for a vertex $v \in Q_0$ is a disjoint union of Dynkin quivers.

\begin{prop}[{\cite[Proposition 2.4]{GL1}, cf.~\cite[Proposition 2.24]{STW}}]\label{prop : derived equivalence}
There exist triangle equivalences 
\begin{equation*}
\D^b (\PP_A^1) \cong \D^b (\CC \widetilde{\TT}_A) \cong \D^b (\CC Q_A),
\end{equation*}
where $Q_A$ is the extended Dynkin quiver given as follows:
\begin{table}[h]
\centering
\begin{tabular}{|c||c|c|c|c|c|} \hline
$A$ & $(1, p, q)$ & $(2, 2, r)$ & $(2, 3, 3)$ & $(2, 3, 4)$ & $(2, 3, 5)$ \\ \hline
$Q_A$ & $A_{p, q}^{(1)}$ & $D_{r + 2}^{(1)}$ & $E_6^{(1)}$ & $E_7^{(1)}$ & $E_8^{(1)}$ \\ \hline
\end{tabular}

\vspace{5pt}
\caption{Extended Dynkin quivers}
\label{table : Extended Dynkin quivers}
\end{table}

\qed
\end{prop}

The action of Proposition~\ref{prop : braid} is transitive due to Meltzer \cite{M1}.
Namely, these three triangulated categories are mutation equivalent.

Since $\coh(\PP_A^1)$ is hereditary (Proposition~\ref{prop : hereditary}), 
indecomposable objects and exceptional objects in $\D^b(\PP_A^1)$ are
indecomposable {\em sheaves} and the exceptional sheaves up to translations.
The Auslander--Reiten translation on $\D^b (\PP_A^1)$ is given by $\S[-1]=- \otimes \O(\vec{\omega})$,
which in particular respects the abelian category $\coh(\PP_A^1)$.
Moreover, the {\em Auslander--Reiten quiver} of $\D^b (\CC Q_A)$ can be described very explicitly.
See {\cite[Section 4 and Section 5]{H}} for the definition of Auslander--Reiten quivers including $\ZZ Q_{A}$.
\begin{prop}[{\cite[Section 5.5]{H}, \cite[Section 5.4.1]{GL1}}]\label{prop : AR quiver}
The components of the Auslander--Reiten quiver are of the form $\ZZ Q_A$ and $\ZZ A_{\infty / r}$ for some positive integer $r$.
In particular, one has
\begin{equation*}
\Ind(\coh(\PP_A^1)) = \left( \ZZ Q_A \sqcup \Big( \bigsqcup_{i = 1}^3 \ZZ A_{\infty / a_i} \Big) \sqcup \Big( \bigsqcup_{\lambda \in \PP^1 \setminus \{ 0, 1, \infty \}} \ZZ A_{\infty / 1} \Big) \right)_0,
\end{equation*}
where $(-)_0$ denotes the set of vertices.
More precisely, the following holds:
\begin{enumerate}
\item The component of the Auslander--Reiten quiver containing vector bundles on $\PP_A^1$ is of type $\ZZ Q_A$.
\item For each $i = 1, 2, 3$, the component of the Auslander--Reiten quiver containing $S_{{\color{red} (i, j)}}$ for $j = 0, \dots, a_i - 1$ is of type $\ZZ A_{\infty / a_i}$.
\item For $\lambda \in \PP^1 \setminus \{ 0, 1, \infty \}$, the component of the Auslander--Reiten quiver containing $S_\lambda$  is of type $\ZZ A_{\infty / 1}$.
\qed
\end{enumerate}
\end{prop}

The exceptional sheaves are now classified rather easily as follows:
\begin{prop}[{cf.~\cite[Section 3.2]{M}}]\label{prop : exceptional object}
If an indecomposable sheaf $E \in \coh(\PP_A^1)$ is exceptional, then $E$ is a vector bundle or a finite length sheaf supported on an orbifold point 
$\lambda_i$ for some $i = 1, 2, 3$.
Moreover, the following holds:
\begin{enumerate}
\item A vector bundle on $\PP_A^1$ is exceptional if and only if it is indecomposable.
\item An indecomposable finite length sheaf $E \in \U_{\lambda_i}$ for some $i = 1, 2, 3$ is exceptional if and only if the length of $E$ is smaller than $a_i$.
\qed
\end{enumerate}
\end{prop}

\section{The number of full exceptional collections}\label{sec : main result}

\subsection{Dynkin case}
Let $\vd=(\Delta_0,\Delta_1)$ be a Dynkin quiver where $\Delta_0$ (resp. $\Delta_1$) denote the set of vertices (resp. arrows), the derived category $\D^b(\CC\vd)$ of finitely generated modules over the path algebra $\CC\vd$. We do not specify an orientation since the derived category of a Dynkin quiver is independent of the choice of orientations.
The Dynkin diagram, the underlying graph of $\vd$, will be denoted by $\Delta$.
The Grothendieck group $K_0(\D^b(\CC\vd))$ of $\D^b(\CC\vd)$ is a free abelian group of rank $|\Delta_0|$, which we shall call the rank of $\vd$ in this paper for simplicity.

For a Dynkin quiver of rank $\mu$, define the number $e(\D^b(\CC\vd)) \in \ZZ$ by 
\begin{equation}\label{eq : the number of full exceptional collections for Dynkin cases}
e(\D^b(\CC\vd)) \coloneqq | \FEC(\D^b(\CC\vd)) / \ZZ^\mu |,
\end{equation}
where the $\ZZ^{\mu}$-action on $\FEC(\D^b(\CC\vd))$ is the one defined in Proposition \ref{prop : braid}.
Since the abelian category $\mod(\CC \vd)$ is hereditary, we have 
\[
e(\D^b(\CC\vd)) =\left| \{(E_1,\dots, E_\mu)\in \FEC(\D^b(\CC\vd))\,|\, E_1,\dots, E_\mu\in \mod(\CC\vd)\} \right|.
\]

Obaid--Nauman--Shammakh--Fakieh--Ringel studied the RHS and obtained the following recursive formula, 
which is originally given by Deligne in a seemingly different but actually the equivalent context.
\begin{prop}[{\cite[Corollaire 1.5]{De}, \cite[Section 4]{ONSFR}}]\label{prop : recursion formula for Dynkin cases}
Let $h$ be a positive number called the Coxeter number, namely, the order of the automorphism $[\S[-1]]\in{\rm Aut}(K_0(\D^b(\CC\vd)))$ called the Coxeter transformation. 
We have 
\begin{equation*}
e(\D^b(\CC\vd)) = \dfrac{h}{2} \sum_{v \in \Delta_0} e(\D^b(\CC\vd^{(v)})),
\end{equation*}
where $\vd^{(v)}$ the full subquiver of $\vd$ given by removing the vertex $v$ and arrows connecting with $v$, 
which is a union of Dynkin quivers of smaller ranks.
\begin{table}[h]
\centering
\begin{tabular}{|c||c|c|c|c|c|} \hline
$\Delta$ & $A_\mu$ & $D_\mu$ & $E_6$ & $E_7$ & $E_8$ \\ \hline
$h$ & $\mu + 1$ & $2(\mu - 1)$ & $12$ & $18$ & $30$ \\ \hline
\end{tabular}

\vspace{5pt}
\label{table : Coxeter number}\caption{Coxeter number}
\end{table}

\qed
\end{prop}
\begin{rem}
Note that $(\S[-1])^h=[2]$, which yields the factor $h/2$ in the formula.
\end{rem}

\begin{prop}[{\cite{De, ONSFR, S-U}}]\label{prop : number of full exceptional collections for Dynkin}
Let $\vd$ be a Dynkin quiver of rank $\mu$. 
We have 
\begin{equation}\label{eq:LL fo ADE}
e(\D^b(\CC\vd)) = \frac{\mu!}{d_1\cdots d_\mu}h^\mu.
\end{equation}
where $2\le d_1\le\dots\le d_\mu=h$ are degrees of $\mu$ algebraic independent invariants of the Weyl group for the Dynkin diagram $\Delta$.
\begin{table}[h]
\centering
\begin{tabular}{|c||c|c|c|c|c|} \hline
$\Delta$ & $A_\mu$ & $D_\mu$ & $E_6$ & $E_7$ & $E_8$ \\ \hline
$e(\D^b(\CC\vd))$ & $(\mu + 1)^{\mu - 1}$ & $2(\mu - 1)^\mu$ & $2^9 \cdot 3^4$ & $2 \cdot 3^{12}$ & $2 \cdot 3^5 \cdot 5^7$ \\ \hline
\end{tabular}

\vspace{5pt}
\label{table : number of full exceptional collections for Dynkin}
\caption{Number of full exceptional collections for Dynkin quivers}
\end{table}

\qed
\end{prop}
In singularity theory, the number in the RHS of \eqref{eq:LL fo ADE} is known as the the degree of the Lyashko--Looijenga map for a simple singularity and 
Deligne shows that it coincides with the number of the set of distinguished bases of vanishing cycles of the singularity modulo signs (\cite{De}, see also~\cite{HR}).
On the other hand, it is proven by \cite{S-P} that the derived Fukaya--Seidel category of a simple singularity is equivalent to 
the derived category $\D^b(\CC\vd)$ of the Dynkin quiver of the corresponding type.
In particular, the equivalence induces a natural bijection between the set of distinguished bases of vanishing cycles of the singularity modulo signs 
and the set $\FEC(\D^b(\CC\vd)) / \ZZ^\mu$.

\subsection{The case of orbifold projective lines}

Obviously, $\FEC(\D^b(\PP_A^1)) / \ZZ^\mu $ is an infinite set. 
It is natural to consider an action of an infinite group on it to obtain a finite number generalizing the definition \eqref{eq : the number of full exceptional collections for Dynkin cases}.
The key idea is the use of spherical twists.

Let $\D \coloneqq \D^b(\PP^1_{A'})\times \D^b(\CC\vd_1)\times \dots \times \D^b(\CC\vd_{k'})$ 
for some $A'$ and Dynkin quivers $\vd_1, \dots, \vd_{k'}$.
\begin{defn}[{\cite[Definition 1.1]{ST}}]
An object $S \in \D$ is called {\em spherical} if $\S(S)\cong S[1]$ where $\S$ is the Serre functor of $\D$ and 
\begin{equation*}
\Hom_{\D} (S, S[p]) \cong 
\begin{cases}
\CC, & p = 0, 1, \\
0, & p \ne 0, 1.
\end{cases}
\end{equation*}
\end{defn}
For the convenience, we denote by $\Sph(\D)$ the set of isomorphism classes of spherical objects in $\D$.

\begin{prop}[{\cite[Proposition 2.10]{ST}}]
For a spherical object $S \in \D$, there exists an autoequivalence ${\rm Tw}_S \in \Aut(\D)$ defined by the exact triangle 
\begin{subequations}
\begin{equation}\label{eq : spherical twist}
\RR \Hom_{\D} (S, E) \otimes S \longrightarrow E \longrightarrow {\rm Tw}_S(E)
\end{equation}
for any object $E \in \D$.
The inverse functor ${\rm Tw}_S^{-1} \in \Aut(\D)$ is given by 
\begin{equation}
{\rm Tw}_S^{-1}(E) \longrightarrow E \longrightarrow S \otimes \RR \Hom_{\D^b(\PP_A^1)} (E, S)^*.
\end{equation}
\end{subequations}
\qed
\end{prop}

We define a subgroup $\ST(\D)$ of the group of autoequivalences $\Aut(\D)$, which stands for the initial of Seidel--Thomas or spherical twists.
\begin{defn}
Define a group $\ST(\D)$ as the set generated by spherical twists:
\begin{equation}
\ST(\D) \coloneqq \big\langle {\rm Tw}_S \in \Aut(\D) \mid S \in \Sph(\D) \big\rangle .
\end{equation}
\end{defn}

From the view point of homological mirror symmetry, the group $\ST(\D)$ is closely related to the mapping class group of a Riemannian surface.
Indeed, we may take the above $\D$ as 
a derived category of the partially wrapped Fukaya category of a marked bordered surface, a derived category of (skew) gentle algebra and so on, 
and consider this group and also the number $e$ defined in \eqref{eq : the number of full exceptional collections for orbifold projective lines} below 
(but we will not discuss on this anymore in this paper).
The {\em Alexander method} is a powerful tool to compute a mapping class group.
As an analogue of this, we prove the following
\begin{lem}[Alexander method]\label{lem : Alexander method}
Let $E \in \D$ be an exceptional object and $S \in \D$ a spherical object.
We have ${\rm Tw}_S (E) \cong E$ if and only if $S \in E^\perp$.
In particular, it holds that
\[
\ST(E^\perp) \cong \big\langle {\rm Tw}_S \in \ST(\D) \mid S \in \Sph(\D), ~ {\rm Tw}_S (E) \cong E \big\rangle .
\]
\end{lem}
\begin{pf}
By the Serre duality, we have 
\[
\Hom_{\D} (S, E) \cong \Hom_{\D} (E, \S(S) ) \cong \Hom_{\D} (E, S[1]).
\]
Hence, $\RR\Hom_{\D} (S, E) = 0$ if and only if $\RR\Hom_{\D} (E, S) = 0$, which yields the statement.
\qed
\end{pf}
Here, due to \cite[Theorem~3.3.2 and Theorem~3.3.3]{M}, $E^\perp$ is of the form 
$\D^b(\PP_{A'}^1)\times \D^b(\CC\vd_1)\times \dots \times \D^b(\CC\vd_{k'})$ for some $A'$ and Dynkin quivers $\vd_1, \dots, \vd_{k'}$.

It is important to note that $\ST(\D^b(\CC\vd)) = \{ 1 \}$ for a Dynkin quiver $\vd$ since $\D^b(\CC\vd)$ has no spherical objects.
Since in this paper we always assume that $\chi_A>0$, we have the following
\begin{prop}\label{prop : Braid group of an orbifold projective line}
We have $\ST(\D^b(\PP_A^1)) = \langle - \otimes \O(\vec{c}) \rangle \cong \ZZ$.
\end{prop}
\begin{pf}
By Proposition \ref{prop : AR quiver} and \ref{prop : exceptional object}, we obtain 
\[
\Sph(\D^b(\PP_A^1)) = \{ S_\lambda [p] \in \D^b(\PP_A^1) \mid \lambda \in \PP^1 \setminus \{ 0, 1, \infty \}, ~ p \in \ZZ \}.
\]
Note that ${\rm Tw}_{S_\lambda} \cong {\rm Tw}_{(S_\lambda [p])}$ for any $p \in \ZZ$.
We shall prove ${\rm Tw}_{S_\lambda} \cong - \otimes \O(\vec{c})$ for all $\lambda \in \PP^1 \setminus \{ 0, 1, \infty \}$.
Since we have $\RR\Hom (S_\lambda, \O(j\vec{x}_i)) \cong \CC[-1]$, it follows from the exact triangle \eqref{eq : simple at an ordinary point} that ${\rm Tw}_{S_\lambda} (\O(j \vec{x}_i)) \cong \O(j \vec{x}_i) \otimes \O(\vec{c})$ for each $i = 1, 2, 3$ and $j = 0, \dots, a_i - 1$.
In particular, these isomorphisms are given explicitly and one can easily check that they provide natural transformation between
the functors on the full strongly exceptional collection.
By Proposition \ref{prop : full strongly exceptional collection}, we obtain ${\rm Tw}_{S_\lambda} \cong - \otimes \O(\vec{c})$.
\qed
\end{pf}

Since the group $\ST(\D^b(\PP_A^1))$ is a subgroup of $\Aut(\D^b(\PP_A^1))$, the set $\FEC(\D^b(\PP_A^1))$ admits the $\ST(\D^b(\PP_A^1))$-action.
Therefore, in view of Lemma~\ref{lem : Alexander method}, we define $e(\D^b(\PP_A^1))$ by 
\begin{equation}\label{eq : the number of full exceptional collections for orbifold projective lines}
e(\D^b(\PP_A^1)) \coloneqq |\FEC(\D^b(\PP_A^1)) / \langle \ST(\D^b(\PP_A^1)) , \ZZ^{\mu_A}\rangle |,
\end{equation}
which is a natural generalization of the definition \eqref{eq : the number of full exceptional collections for Dynkin cases}.

Recall that the abelian category $\coh(\PP_A^1)$ is hereditary (Proposition \ref{prop : hereditary}) and the group
$\ST(\D^b(\PP_A^1))$ respects $\coh(\PP_A^1)$ due to Proposition~\ref{prop : Braid group of an orbifold projective line}, we have 
\[
e(\D^b(\PP_A^1)) =\left| \{(E_1,\dots, E_\mu)\in \FEC(\D^b(\PP_A^1))\,|\, E_1,\dots, E_\mu\in \coh(\PP_A^1)\} \left/\langle - \otimes \O(\vec{c})\rangle\right.\right|.
\]

Lemma~\ref{lem : Alexander method} enable us to have the following recursive formula as an analogue of Proposition \ref{prop : recursion formula for Dynkin cases}.
\begin{thm}\label{thm : recursion formula}
We have
\begin{equation}\label{eq : recursive formula}
e(\D^b(\PP_A^1)) = \dfrac{1}{\chi_A} \sum_{v \in (Q_A)_0} e(\D^b(\CC Q_A^{(v)})) + \sum_{i = 1}^3 a_i \sum_{j = 1}^{a_i - 1} 
\binom{\mu_A-1}{a_i - j-1}\cdot
e(\D^b(\PP^1_{A_{(i, j)}}) )\cdot e(\D^b(\CC\vec{A}_{a_i - j - 1})),
\end{equation}
where $A_{(i, j)} = (a'_1, a'_2, a'_3)$ is defined by $a'_k: = a_k$ if $k\ne i$ and $a'_k  \coloneqq  j$ if $k=i$.
\end{thm}
Hence, we see from this recursive formula that $e(\D^b(\PP_A^1)) $ is indeed a positive integer.
The proof of Theorem~\ref{thm : recursion formula} will be given in Section \ref{sec : proof of main thm 1}.

\begin{thm}\label{thm : main}
Let $\PP_A^1$ be an orbifold projective line with $\chi_A>0$.
We have 
\begin{equation}\label{eq:LL affine}
e(\D^b(\PP_A^1)) = \frac{\mu_A!}{a_1!a_2!a_3!\chi_A}a_1^{a_1}a_2^{a_2}a_3^{a_3} .
\end{equation}
\end{thm}
The proof of Theorem~\ref{thm : main} will be given in Section \ref{sec : proof of main thm 2}.

Bridgeland introduced the notion of a stability condition on a triangulated category \cite{Br}.
He showed that the space $\Stab(\D)$ of all stability conditions on a triangulated category $\D$ is a complex manifold.
We expect that the quotient space $\Stab(\D^b(\PP_A^1)) / \ST(\D^b(\PP_A^1))$ is a complex manifold isomorphic to $M=\CC^{\mu_A-1}\times\CC^*$ 
on which the Frobenius structure from the Gromov--Witten theory for $\PP_A^1$ lives. 
On the other hand, via the mirror symmetry between $\PP^{1}_{A}$ and the corresponding affine cusp polynomial $f_{A}$
explained in Section~\ref{LL-Frob}, the quotient space $\Stab(\D^b(\PP_A^1)) / \ST(\D^b(\PP_A^1))$ is expected to be mirror isomorphic to the parameter space of the universal unfolding for $f_{A}$ on which the Frobenius structure from the Saito theory for $f_{A}$ also lives.

There, elements in $\FEC(\D^b(\PP_A^1))$ are viewed as distinguished bases of vanishing cycles 
as similar to cycles for isolated hypersurface singularities 
and the braid group actions on $\FEC(\D^b(\PP_A^1))$ are regarded as its geometric actions on distinguished bases (see \cite{E}).
Note that pairwise distinct critical values play an important role to make such distinguished bases and have the interpretation 
as the eigenvalues of the multiplication operator by the Euler vector field for the Frobenius manifold.
This is the reason why an open dense set where critical values are pairwise distinct in the parameter space, 
the complement of the bifurcation set, relates to our purely categorical quantity.

The braid group acts on $\FEC(\D^b(\PP_A^1)) / \langle \ST(\D^b(\PP_A^1)), \ZZ^{\mu_A} \rangle$
since these actions commute with autoequivalence group actions on $\FEC(\D^b(\PP_A^1))$.
Namely, the set $\FEC(\D^b(\PP_A^1)) / \langle \ST(\D^b(\PP_A^1)), \ZZ^{\mu_A} \rangle$ is
represented by the quotient set of the braid group by the stabilizer group for an element in 
$\FEC(\D^b(\PP_A^1)) / \langle \ST(\D^b(\PP_A^1)), \ZZ^{\mu_A} \rangle$ since the action of braid group is transitive
as explained after Proposition \ref{prop : derived equivalence}.
It is also expected that this stabilizer group should coincide with the fundamental group for the complement of the bifurcation set. 
Indeed, the number in the RHS of \eqref{eq:LL affine} is known as the degree of 
the Lyashko--Looijenga map for the Frobenius structure (see Section~\ref{sec : generalized LL map} and \cite{DZ, He}).
Theorem~\ref{thm : main} hints a consistency among mirror symmetry, space of stability conditions and Duvrovin's conjecture.

\section{Proof of Theorem \ref{thm : recursion formula}}\label{sec : proof of main thm 1}

For simplicity, we write $\Exc( \coh (\PP_A^1)) / \O(\vec{c})$ to mean the quotient set $\Exc( \coh (\PP_A^1)) /\langle(-)\otimes \O(\vec{c})\rangle$, and similar for others.
By the definition of $e(\D^b(\PP_A^1))$, it follows from Lemma \ref{lem : Alexander method} and Proposition \ref{prop : Braid group of an orbifold projective line} that
\begin{equation*}
e (\D^b(\PP_A^1)) = \sum_{E \in \Exc( \coh (\PP_A^1)) / \O(\vec{c})} e (E^\perp),
\end{equation*}
where $e (E^\perp)$ is defined as 
\begin{equation*}
e(E^\perp) \coloneqq |\FEC(E^\perp) / \langle \ST(E^\perp) , \ZZ^{\mu_A-1}\rangle |,
\end{equation*}

Hence, Proposition \ref{prop : exceptional object} implies that 
\begin{equation}\label{eq : formula1}
\sum_{E \in \Exc( \coh (\PP_A^1)) / \O(\vec{c})} e (E^\perp) = \sum_{E \in \Exc( {\rm vect} (\PP_A^1)) / \O(\vec{c})} e (E^\perp) + \sum_{i = 1}^3 \sum_{E \in \Exc( \U_{\lambda_i}) / \O(\vec{c})} e (E^\perp).
\end{equation}
\begin{lem}\label{lem : classification of exceptional objects}
There exist one-to-one correspondences
\begin{subequations}
\begin{equation}\label{eq : exceptional vector bundles}
(Q_A)_0 \xrightarrow{1 \colon 1} \Exc({\rm vect}(\PP_A^1)) / \O(\vec{\omega}), \quad v \mapsto [E^{Q_A}_v]
\end{equation}
\begin{equation}\label{eq : finite length exceptional sheaves}
(\widetilde{\TT}_A)_0 \setminus \{ \mathbf{1}, \mathbf{1}^* \} \xrightarrow{1 \colon 1} \Exc(\coh_0 (\PP_A^1)) / \O(\vec{\omega}), \quad (i, j) \mapsto [E^{\widetilde{\TT}_A}_{(i, j)}].
\end{equation}
\end{subequations}
\end{lem}
\begin{pf}
First, we prove \eqref{eq : exceptional vector bundles}.
It was proved by \cite[Section 5.5]{H} and \cite[Section 5.4.1]{GL1} that the component of the Auslander--Reiten quiver containing vector bundles on $\PP_A^1$ is of type $\ZZ Q_A$.
By the definition of the Auslander--Reiten quiver, the statement follows from (i) of Proposition \ref{prop : exceptional object}.

Next, we prove \eqref{eq : finite length exceptional sheaves}.
For each $i$, we can choose the equivalence $\D^b (\PP_A^1) \cong \D^b (\CC \widetilde{\TT}_A)$ so that the simple module $S_{i,a_i-1}$ is mapped to 
$E^{\widetilde{\TT}_A}_{(i, a_i-1)}$ and $E^{\widetilde{\TT}_A}_{(i, 1)}\longrightarrow \dots \longrightarrow E^{\widetilde{\TT}_A}_{(i, a_i-1)}$ 
belong to the component of the Auslander--Reiten quiver of type $\ZZ A_{\infty / a_i}$. 
Here note that $E^{\widetilde{\TT}_A}_{(i, j)}$ are obtained by iterated extensions of $S_{(i,j)}$s, in particular,
$E^{\widetilde{\TT}_A}_{(i, j)}$ is not isomorphic to $S_{(i,j)}$ except for $j=a_{i}-1$ since ${\rm Hom} (S_{(i,j)}, S_{(i,j')})=0$ if $j\neq j'$.
Since $S_{i,j}\otimes\O(\vec{\omega})\cong S_{i,j-1}$, (ii) of Proposition \ref{prop : exceptional object} implies that the number of the set $\Exc(\coh_0 (\PP_A^1)) / \O(\vec{\omega})$ is equal to the number of the set $(\widetilde{\TT}_A)_0 \setminus \{ \mathbf{1}, \mathbf{1}^* \} $.
Hence, the correspondence is bijective.
\qed
\end{pf}

\begin{prop}\label{prop : perpendicular category}
Let $E \in \D^b(\PP_A^1)$ be an exceptional object.
\begin{enumerate}
\item For each vertex $v \in (Q_A)_0$, we have 
\begin{equation*}
(E^{Q_A}_v)^\perp \cong \D^b(\CC Q_A^{(v)}).
\end{equation*}
\item For each vertex $(i, j) \in (\widetilde{\TT}_A)_0 \setminus \{ \mathbf{1}, \mathbf{1}^* \}$, we have
\begin{equation*}
(E^{\widetilde{\TT}_A}_{(i, j)})^\perp \cong \D^b(\PP^1_{A_{(i, j)}}) \times \D^b(\CC\vec{A}_{a_i - j - 1}).
\end{equation*}
\end{enumerate}
\end{prop}
\begin{pf}
We embed the quiver $Q_A$ into the component of the Auslander--Reiten quiver of type $\ZZ Q_A$ so that the vertex $v$ is mapped to $E^{Q_A}_v$.
The first statement follows from the fact that we may choose a subquiver $Q'$ of $\ZZ Q_A$ giving rise to another full exceptional collection in $\D^b(\CC Q_A)$, 
by ``bending suitable arrows", whose underlying graph is isomorphic to the one of $Q_A$ in which the vertex $v$ is a sink corresponding to the object $E^{Q_A}_v$.

The second statement follows from that we may replace in the component of the Auslander--Reiten quiver of type $\ZZ A_{\infty / a_i}$ 
the subquiver $E^{\widetilde{\TT}_A}_{(i, 1)}\longrightarrow \dots \longrightarrow E^{\widetilde{\TT}_A}_{(i, j)}\longrightarrow \dots \longrightarrow E^{\widetilde{\TT}_A}_{(i, a_i-1)}=S_{i,a_i-1}$ with the subquiver $E^{\widetilde{\TT}_A}_{(i, 1)}\longrightarrow \dots \longrightarrow E^{\widetilde{\TT}_A}_{(i, j)}\longleftarrow \dots \longleftarrow E^{\widetilde{\TT}_A}_{(i, a_i-1)}=S_{i,a_i-1}$ 
to get another full exceptional collection in $\D^b (\CC \widetilde{\TT}_A)$. 
\qed
\end{pf}
\begin{lem}
We have 
\begin{equation}\label{eq : counting of exceptional vector bundle}
\sum_{E \in \Exc( {\rm vect} (\PP_A^1)) / \O(\vec{c})} e (E^\perp) = \dfrac{1}{\chi_A} \sum_{v \in Q_A} e(\D^b(\CC Q_A^{(v)})).
\end{equation}
\end{lem}
\begin{pf}
Let $a$ be the least common multiple of $a_1, a_2$ and $a_3$. Since $a\vec{\omega} = - a\chi_A \vec{c}$ in $L_A$, we have 
\[
\sum_{E \in \Exc( {\rm vect} (\PP_A^1)) / \O(a \chi_A \vec{c})} e (E^\perp) = a \chi_A \sum_{E \in \Exc( {\rm vect} (\PP_A^1)) / \O(\vec{c})} e (E^\perp)
\]
and
\[
\sum_{E \in \Exc( {\rm vect} (\PP_A^1)) / \O(a \vec{\omega})} e (E^\perp) = a \sum_{E \in \Exc( {\rm vect} (\PP_A^1)) / \O(\vec{\omega})} e (E^\perp)
\]
Hence, we obtain 
\[
\sum_{E \in \Exc( {\rm vect} (\PP_A^1)) / \O(\vec{c})} e (E^\perp) = \dfrac{1}{\chi_A} \sum_{E \in \Exc( {\rm vect} (\PP_A^1)) / \O(\vec{\omega})} e (E^\perp).
\]
Lemma \ref{lem : classification of exceptional objects} and Proposition \ref{prop : perpendicular category} imply the statement.
\qed
\end{pf}

The following lemma computes the number of full exceptional collections when $E^\perp $ is a product of two triangulated categories.
\begin{lem}[Shuffle Lemma]\label{lem : shuffle lemma}
Let $Q_1$ and $Q_2$ be Dynkin quivers or extended Dynkin quivers.
Denote by $\mu_i$ the length of a full exceptional collection in $\D^b(\CC Q_i)$. 
We have
\[
e(\D^b(\CC Q_1) \times \D^b(\CC Q_2)) = \dfrac{(\mu_1 + \mu_2)!}{\mu_1! \, \mu_2!}\cdot e(\D^b(\CC Q_1)) \cdot e(\D^b(\CC Q_2)).
\]
\end{lem}
\begin{pf}
We can show the statement by direct calculations (cf.~\cite[Section 2]{ONSFR}).
\qed
\end{pf}

\begin{lem}
For every $i = 1, 2, 3$, we have 
\begin{equation}\label{eq : counting of exceptional finite length sheaf}
\sum_{E \in \Exc( \U_{\lambda_i}) / \O(\vec{c})} e (E^\perp) = a_i \sum_{j = 1}^{a_i - 1} \binom{\mu_A-1}{a_i-j-1}\cdot  e(\D^b(\PP^1_{A_{(i, j)}}))\cdot e(\D^b(\CC \vec{A}_{a_i - j - 1})).
\end{equation}
\end{lem}
\begin{pf}
The autoequivalence $- \otimes \O(\vec{c})$ acts trivially on $\Exc(\U_{\lambda_i})$.
On the other hand, the order of the autoequivalence $- \otimes \O(\vec{\omega})$ acting on $\Exc(\U_{\lambda_i})$ is $a_i$.
Hence, there is an $a_i$-to-one correspondence between $\Exc(\U_{\lambda_i}) / \O(\vec{c})$ and $\Exc( \U_{\lambda_i}) / \O(\vec{\omega})$, thus
\[
\sum_{E \in \Exc( \U_{\lambda_i}) / \O(\vec{c})} e (E^\perp) = a_i \sum_{E \in \Exc( \U_{\lambda_i}) / \O(\vec{\omega})} e (E^\perp).
\]
Therefore, Lemma~\ref{lem : classification of exceptional objects}, Proposition~\ref{prop : perpendicular category} and Lemma~\ref{lem : shuffle lemma} imply the statement.
\qed
\end{pf}

Putting the equations \eqref{eq : formula1}, \eqref{eq : counting of exceptional vector bundle} and \eqref{eq : counting of exceptional finite length sheaf} together, we obtain the equation \eqref{eq : recursive formula}.
\qed

\section{Proof of Theorem \ref{thm : main}}\label{sec : proof of main thm 2}

We prove the statement by Theorem \ref{thm : recursion formula} and induction, on a case-by-case basis.

\subsection{Case of $A = (1, p, q)$}
We prove the statement by induction on $\mu_A=p+q$.
When $\mu_A = 1$, namely $A = (1, 1, 1)$, we have $e(\D^b(\PP^1)) = 1$.
Hence, we assume that $\mu_A \ge 2$.
For each vertex $v \in (A_{p, q}^{(1)})_0$, the quiver $(A_{p, q}^{(1)})^{(v)}$ is isomorphic to a Dynkin quiver of type $A_{p+q-1}$.
\begin{lem}[{\cite{Hu}, cf.~\cite[Theorem 3]{S-V}}]\label{lem : Hurwitz1}
For positive numbers $p,q$, we have the following identity:
\begin{align*}
\begin{split}
& \dfrac{(p + q - 1)!}{(p - 1)! (q - 1)!} p^p q^q \\
= ~ & ~ p q (p + q)^{p + q - 2} \\
& + p \sum_{j = 1}^{p - 1} \frac{(p + q - 1)!}{(q + j)! (p - j - 1)!} \frac{(q + j - 1)!}{(j - 1)! (q - 1)!} j^j q^q (p - j)^{p - j - 2} \\
& + q \sum_{j = 1}^{q - 1} \frac{(p + q - 1)!}{(p + j)! (q - j - 1)!} \frac{(p + j - 1)!}{(p - 1)! (j - 1)!} p^p j^j (q - j)^{q - j - 2} .
\end{split}
\end{align*}
\qed
\end{lem}
Since the Euler character is given by $\chi_A = (p + q) / p q$, Theorem~\ref{thm : recursion formula} and Lemma~\ref{lem : Hurwitz1} 
imply
\begin{eqnarray*}
e(\D^b(\PP_{1,p,q}^1)) & = & p q (p + q)^{p + q - 2} \\
& & + p \sum_{j = 1}^{p - 1} \frac{(p + q - 1)!}{(q + j)! (p - j - 1)!} \frac{(q + j - 1)!}{(j - 1)! (q - 1)!} j^j q^q (p - j)^{p - j - 2} \\
& & + q \sum_{j = 1}^{q - 1} \frac{(p + q - 1)!}{(p + j)! (q - j - 1)!} \frac{(p + j - 1)!}{(p - 1)! (j - 1)!} p^p j^j (q - j)^{q - j - 2} \\
& = & \frac{(p + q - 1)!}{(p - 1)! (q - 1)!} p^p q^q \\
& = & \frac{(p + q)!}{p! q! \frac{p + q}{pq}} p^p q^q.
\end{eqnarray*}
Hence, we obtain the statement for $A = (1, p, q)$.

\subsection{Case of $A = (2, 2, r)$}
We prove the statement by induction on $r$.
When $r = 1$, the statement is already proven. Assume that $r \ge 2$.
Here, we use Proposition \ref{prop : number of full exceptional collections for Dynkin} and Theorem \ref{thm : main} for $A = (1, p, q)$ to compute these numbers.
\begin{table}[hbtp]
\centering
{\small
\begin{tabular}{|c|c|c|} \hline
$v \in (Q_A)_0$ & $\D^b(\CC Q_A^{(v)})$ & $e(\D^b(\CC Q_A^{(v)}))$ \\ \hline \hline
$v = 1$ & $\D^b(\CC \vec{D}_{r + 2})$ & $2 (r + 1)^{r + 2}$ \\ \hline
$v = 2$ & $\D^b(\CC \vec{D}_{r + 2})$ & $2 (r + 1)^{r + 2}$ \\ \hline
$v = 3$ & $\D^b(\CC \vec{D}_{2}) \times \D^b(\CC \vec{D}_{r})$ & $\dfrac{(r + 2)!}{2! r!} 2^2 (r - 1)^r$ \\ \hline
$\vdots$ & $\vdots$ & $\vdots$ \\ \hline
$v = k$ & $\D^b(\CC \vec{D}_{k - 1}) \times \D^b(\CC \vec{D}_{r + 3 - k})$ & $\dfrac{(r + 2)!}{(k - 1)! (r + 3 - k)!} 2^2 (k - 2)^{k - 1} (r + 2 - k)^{r + 3 - k}$ \\ \hline
$\vdots$ & $\vdots$ & $\vdots$ \\ \hline
$v = r + 1$ & $\D^b(\CC \vec{D}_{r}) \times \D^b(\CC \vec{D}_{2})$ & $\dfrac{(r + 2)!}{r! 2!} 2^2 (r - 1)^r$ \\ \hline
$v = r + 2$ & $\D^b(\CC \vec{D}_{r + 2})$ & $2 (r + 1)^{r + 2}$ \\ \hline
$v = r + 3$ & $\D^b(\CC \vec{D}_{r + 2})$ & $2 (r + 1)^{r + 2}$ \\ \hline
\end{tabular}

\vspace{10pt}
\centering
\begin{tabular}{|c|c|c|} \hline
$v \in (\widetilde{\TT}_A)_0$ & $\D^b(\PP_{A_{(i,j)}}^1) \times \D^b(\CC \vec{A}_{a_i - j - 1})$ & $\binom{\mu_A-1}{a_i-j-1}\cdot  e(\D^b(\PP^1_{A_{(i, j)}}))\cdot e(\D^b(\CC \vec{A}_{a_i - j - 1}))$ \\ \hline \hline
$v = (1,1)$ & $\D^b(\PP^1_{1, 2, r})$ & $4 (r + 1) r^{r + 1}$ \\ \hline
$v = (2,1)$ & $\D^b(\PP^1_{2, 1, r})$ & $4 (r + 1) r^{r + 1}$ \\ \hline
$v = (3,1)$ & $\D^b(\PP^1_{2, 2, 1}) \times \D^b(\CC \vec{A}_{r - 2})$ & $\dfrac{(r + 2)!}{4! (r - 2)!} 4 \cdot 2 \cdot 3 \cdot 4 \cdot (r - 1)^{r - 3}$ \\ \hline
$\vdots$ & $\vdots$ & $\vdots$ \\ \hline
$v = (3,j)$ & $\D^b(\PP^1_{2, 2, j}) \times \D^b(\CC \vec{A}_{r - j - 1})$ & $\dfrac{(r + 2)!}{(j + 3)! (r - j - 1)!} 4(j + 1)(j + 2)(j + 3)j^{j + 1} (r - j)^{r - j - 2}$ \\ \hline
$\vdots$ & $\vdots$ & $\vdots$ \\ \hline
$v = (3,r - 1)$ & $\D^b(\PP^1_{2, 2, r - 1})$ & $4r(r + 1)(r + 2)(r - 1)^{r}$ \\ \hline
\end{tabular}
}

\vspace{5pt}
\caption{Case of $A = (2, 2, r)$}
\label{table : (2, 2, r)}
\end{table}

\begin{lem}[{\cite{Hu}, cf.~\cite[Theorem 3]{S-V}}]\label{lem : Hurwitz2}
For a positive number $r$, we have the following identity:
\begin{align*}
 \begin{split}
& (r + 1)(r + 2)(r + 3)r^{r + 1} \\
= ~ & ~ 4(r + 1)r^{r + 1} + 2(r + 1)^{r + 2} \\
& + r \sum_{j = 1}^{r - 1} \dfrac{(r + 2)!}{(j + 3)! (r - j - 1)!} 4(j + 1)(j + 2)(j + 3)j^{j + 1} (r - j)^{r - j - 2} \\
& + r \sum_{k = 1}^{r - 1} \frac{(r + 2)!}{(k + 1)! (r - k + 1)!} k^{k + 1} (r - k)^{r - k + 1} .
\end{split}
\end{align*}
\qed
\end{lem}

Since $\chi_A = 1 / r$, Theorem~\ref{thm : recursion formula}, Table \ref{table : (2, 2, r)} and Lemma~\ref{lem : Hurwitz2} imply
\begin{eqnarray*}
e(\D^b(\PP^1_{2,2,r})) & = & r \Big( 8 (r + 1)^{r + 2} + \sum_{k = 3}^{r + 1} \dfrac{(r + 2)!}{(k - 1)! (r + 3 - k)!} 4 (k - 2)^{k - 1} (r + 2 - k)^{r + 3 - k} \Big) \\ 
& & + 2 \cdot 4 (r + 1) r^{r + 1} + 2 \cdot 4 (r + 1) r^{r + 1} \\
& & + r \sum_{j = 1}^{r - 1} \dfrac{(r + 2)!}{(j + 3)! (r - j - 1)!} 4(j + 1)(j + 2)(j + 3)j^{j + 1} (r - j)^{r - j - 2} \\
& = & 16 (r + 1)r^{r + 1} + 8 (r + 1)^{r + 2} \\
& & + 4 r \sum_{j = 1}^{r - 1} \frac{(r + 2)!}{(j - 1)! (r - j + 3)!} (r - j + 1)(r - j + 2)(r - j +3) j^{j - 2} (r - j)^{r - j + 1} \\
& & + 4 r \sum_{k = 1}^{r - 1} \frac{(r + 2)!}{(k + 1)! (r - k + 1)!} k^{k + 1} (r - k)^{r - k + 1} \\
& = & 4 (r + 1) (r + 2) ( r+ 3) r^{r + 1} \\
& = & \frac{(r + 3)!}{2! 2! r! \frac{1}{r}} 2^2 2^2 r^r.
\end{eqnarray*}
Hence, we obtain the statement for $A = (2, 2, r)$.

\subsection{Case of $A = (2, 3, 3)$}
Here, we use Proposition \ref{prop : number of full exceptional collections for Dynkin} and Theorem \ref{thm : main} for $A = (1, p, q)$ and $(2, 2, r)$ to compute these numbers.
\begin{table}[hbtp]
\centering
\begin{tabular}{|c|c|c|} \hline
$v \in (Q_A)_0$ & $\D^b(\CC Q_A^{(v)})$ & $e(\D^b(\CC Q_A^{(v)}))$ \\ \hline \hline
$v = 1$ & $\D^b(\CC \vec{E}_{6})$ & $2^9 \cdot 3^4$ \\ \hline
$v = 2$ & $\D^b(\CC \vec{A}_{1}) \times \D^b(\CC \vec{A}_{5})$ & $6 \cdot 6^4$ \\ \hline
$v = 3$ & $\D^b(\CC \vec{E}_{6})$ & $2^9 \cdot 3^4$ \\ \hline
$v = 4$ & $\D^b(\CC \vec{A}_{1}) \times \D^b(\CC \vec{A}_{5})$ & $6 \cdot 6^4$ \\ \hline
$v = 5$ & $\D^b(\CC \vec{A}_{2}) \times \D^b(\CC \vec{A}_{2}) \times \D^b(\CC \vec{A}_{2})$ & $90 \cdot 3 \cdot 3 \cdot 3$ \\ \hline
$v = 6$ & $\D^b(\CC \vec{A}_{5}) \times \D^b(\CC \vec{A}_{1})$ & $6 \cdot 6^4$ \\ \hline
$v = 7$ & $\D^b(\CC \vec{E}_{6})$ & $2^9 \cdot 3^4$ \\ \hline
\end{tabular}

\vspace{10pt}
\centering
\begin{tabular}{|c|c|c|} \hline
$v \in (\widetilde{\TT}_A)_0$ & $\D^b(\PP_{A_{(i,j)}}^1) \times \D^b(\CC \vec{A}_{a_i - j - 1})$ & $\binom{\mu_A-1}{a_i-j-1}\cdot  e(\D^b(\PP^1_{A_{(i, j)}}))\cdot e(\D^b(\CC \vec{A}_{a_i - j - 1}))$ \\ \hline \hline
$v = (1,1)$ & $\D^b(\PP^1_{1, 3, 3})$ & $21870$ \\ \hline
$v = (2,1)$ & $\D^b(\PP^1_{2, 1, 3}) \times \D^b(\CC \vec{A}_{1})$ & $6 \cdot 1296$ \\ \hline
$v = (2,2)$ & $\D^b(\PP^1_{2, 2, 3})$ & $38880$ \\ \hline
$v = (3,1)$ & $\D^b(\PP^1_{2, 3, 1}) \times \D^b(\CC \vec{A}_{1})$ & $6 \cdot 1296$ \\ \hline
$v = (3,2)$ & $\D^b(\PP^1_{2, 3, 2})$ & $38880$ \\ \hline
\end{tabular}

\vspace{5pt}
\caption{Case of $A = (2, 3, 3)$}
\label{table : (2, 3, 3)}
\end{table}

Since $\chi_A = 1 / 6$, Theorem \ref{thm : recursion formula} and Table \ref{table : (2, 3, 3)} imply 
\begin{eqnarray*}
e(\D^b(\PP^1_A)) & = & 6 (2^9 \cdot 3^4 + 6 \cdot 6^4 + 2^9 \cdot 3^4 + + 6 \cdot 6^4 + 90 \cdot 3^3 + 6 \cdot 6^4 + 2^9 \cdot 3^4) \\ 
& & + 2 \cdot 21870 + 3 (38880 + 6 \cdot 1296) + 3 (38880 + 6 \cdot 1296) \\
& = & 1224720 \\
& = & \frac{7!}{2! 3! 3! \frac{1}{6}} 2^2 3^3 3^3.
\end{eqnarray*}
Hence, we obtain the statement for $A = (2, 3, 3)$.

\subsection{Case of $A = (2, 3, 4)$}
Here, we use Proposition \ref{prop : number of full exceptional collections for Dynkin} and Theorem \ref{thm : main} for $A = (1, p, q)$, $(2, 2, r)$ and $(2, 3, 3)$ to compute these numbers.
\begin{table}[hbtp]
\centering
\begin{tabular}{|c|c|c|} \hline
$v \in (Q_A)_0$ & $\D^b(\CC Q_A^{(v)})$ & $e(\D^b(\CC Q_A^{(v)}))$ \\ \hline \hline
$v = 1$ & $\D^b(\CC \vec{A}_{7})$ & $8^6$ \\ \hline
$v = 2$ & $\D^b(\CC \vec{E}_{7})$ & $2 \cdot 3^{12}$ \\ \hline
$v = 3$ & $\D^b(\CC \vec{A}_{1}) \times \D^b(\CC \vec{D}_{6})$ & $7 \cdot 2 \cdot 5^6$ \\ \hline
$v = 4$ & $\D^b(\CC \vec{A}_{2}) \times \D^b(\CC \vec{A}_{5})$ & $21 \cdot 3 \cdot 6^4$ \\ \hline
$v = 5$ & $\D^b(\CC \vec{A}_{1}) \times \D^b(\CC \vec{A}_{3}) \times \D^b(\CC \vec{A}_{3})$ & $140 \cdot 4^2 \cdot 4^2$ \\ \hline
$v = 6$ & $\D^b(\CC \vec{A}_{5}) \times \D^b(\CC \vec{A}_{2})$ & $21 \cdot 6^4 \cdot 3$ \\ \hline
$v = 7$ & $\D^b(\CC \vec{D}_{6}) \times \D^b(\CC \vec{A}_{1})$ & $7 \cdot 5^6 \cdot 2$ \\ \hline
$v = 8$ & $\D^b(\CC \vec{E}_{7})$ & $2 \cdot 3^{12}$ \\ \hline
\end{tabular}

\vspace{10pt}
\centering
\begin{tabular}{|c|c|c|} \hline
$v \in (\widetilde{\TT}_A)_0$ & $\D^b(\PP_{A_{(i,j)}}^1) \times \D^b(\CC \vec{A}_{a_i - j - 1})$ & $\binom{\mu_A-1}{a_i-j-1}\cdot  e(\D^b(\PP^1_{A_{(i, j)}}))\cdot e(\D^b(\CC \vec{A}_{a_i - j - 1}))$ \\ \hline \hline
$v = (1,1)$ & $\D^b(\PP^1_{1, 3, 4})$ & $60 \cdot 3^3 \cdot 4^4$ \\ \hline
$v = (2,1)$ & $\D^b(\PP^1_{2, 1, 4}) \times \D^b(\CC \vec{A}_{1})$ & $7 \cdot 20 \cdot 2^2 \cdot 4^4$ \\ \hline
$v = (2,2)$ & $\D^b(\PP^1_{2, 2, 4})$ & $860160$ \\ \hline
$v = (3,1)$ & $\D^b(\PP^1_{2, 3, 1}) \times \D^b(\CC \vec{A}_{2})$ & $21 \cdot 1296 \cdot 3$ \\ \hline
$v = (3,2)$ & $\D^b(\PP^1_{2, 3, 2}) \times \D^b(\CC \vec{A}_{1})$ & $7 \cdot 38880$ \\ \hline
$v = (3,3)$ & $\D^b(\PP^1_{2, 3, 3})$ & $1224720$ \\ \hline
\end{tabular}

\vspace{5pt}
\caption{Case of $A = (2, 3, 4)$}
\label{table : (2, 3, 4)}
\end{table}

Since $\chi_A = 1 / 12$, Theorem \ref{thm : recursion formula} and Table \ref{table : (2, 3, 4)} imply
\begin{eqnarray*}
e(\D^b(\PP^1_A)) & = & 12 (8^6 + 2 \cdot 3^{12} + 7 \cdot 2 \cdot 5^6 + 21 \cdot 3 \cdot 6^4 + 140 \cdot 4^2 \cdot 4^2 + 21 \cdot 6^4 \cdot 3 + 7 \cdot 2 \cdot 5^6 + 2 \cdot 3^{12}) \\ 
& & + 2 \cdot 610 \cdot 3^3 \cdot 4^4 + 3 (7 \cdot 20 \cdot 2^2 \cdot 4^4 + 860160) + 4 (21 \cdot 1296 \cdot 3 + 7 \cdot 38840 + 1224720) \\
& = & 46448640 \\
& = & \frac{8!}{2! 3! 4! \frac{1}{12}} 2^2 3^3 4^4.
\end{eqnarray*}
Hence, we obtain the statement for $A = (2, 3, 4)$.

\subsection{Case of $A = (2, 3, 5)$}
Here, we use Proposition \ref{prop : number of full exceptional collections for Dynkin} and Theorem \ref{thm : main} for $A = (1, p, q)$, $(2, 2, r)$, $(2, 3, 3)$ and $(2, 3, 4)$ to compute these numbers.
\begin{table}[hbtp]
\centering
\begin{tabular}{|c|c|c|} \hline
$v \in (Q_A)_0$ & $\D^b(\CC Q_A^{(v)})$ & $e(\D^b(\CC Q_A^{(v)}))$ \\ \hline \hline
$v = 1$ & $\D^b(\CC \vec{A}_{8})$ & $9^7$ \\ \hline
$v = 2$ & $\D^b(\CC \vec{D}_{8})$ & $2 \cdot 7^8$ \\ \hline
$v = 3$ & $\D^b(\CC \vec{A}_{7}) \times \D^b(\CC \vec{A}_{1})$ & $8 \cdot 8^6$ \\ \hline
$v = 4$ & $\D^b(\CC \vec{A}_{1}) \times \D^b(\CC \vec{A}_{2}) \times \D^b(\CC \vec{A}_{5})$ & $168 \cdot 3 \cdot 6^4$ \\ \hline
$v = 5$ & $\D^b(\CC \vec{A}_{4}) \times \D^b(\CC \vec{A}_{4})$ & $70 \cdot 5^3 \cdot 5^3$ \\ \hline
$v = 6$ & $\D^b(\CC \vec{D}_{5}) \times \D^b(\CC \vec{A}_{3})$ & $56 \cdot 2 \cdot 4^5 \cdot 4^2$ \\ \hline
$v = 7$ & $\D^b(\CC \vec{E}_{6}) \times \D^b(\CC \vec{A}_{2})$ & $28 \cdot 2^9 \cdot 3^4 \cdot 3$ \\ \hline
$v = 8$ & $\D^b(\CC \vec{E}_{7}) \times \D^b(\CC \vec{A}_{1})$ & $8 \cdot 2 \cdot 3^{12}$ \\ \hline
$v = 9$ & $\D^b(\CC \vec{E}_{8})$ & $2 \cdot 3^5 \cdot 5^7$ \\ \hline
\end{tabular}

\vspace{10pt}
\centering
\begin{tabular}{|c|c|c|} \hline
$v \in (\widetilde{\TT}_A)_0$ & $\D^b(\PP_{A_{(i,j)}}^1) \times \D^b(\CC \vec{A}_{a_i - j - 1})$ & $\binom{\mu_A-1}{a_i-j-1}\cdot  e(\D^b(\PP^1_{A_{(i, j)}}))\cdot e(\D^b(\CC \vec{A}_{a_i - j - 1}))$ \\ \hline \hline
$v = (1,1)$ & $\D^b(\PP^1_{1, 3, 5})$ & $105 \cdot 3^3 \cdot 5^5$ \\ \hline
$v = (2,1)$ & $\D^b(\PP^1_{2, 1, 5}) \times \D^b(\CC \vec{A}_{1})$ & $8 \cdot 30 \cdot 2^2 \cdot 5^5$ \\ \hline
$v = (2,2)$ & $\D^b(\PP^1_{2, 2, 5})$ & $21000000$ \\ \hline
$v = (3,1)$ & $\D^b(\PP^1_{2, 3, 1}) \times \D^b(\CC \vec{A}_{3})$ & $56 \cdot 12 \cdot 2^2 \cdot 3^3 \cdot 4^2$ \\ \hline
$v = (3,2)$ & $\D^b(\PP^1_{2, 3, 2}) \times \D^b(\CC \vec{A}_{2})$ & $28 \cdot 38880 \cdot 3$ \\ \hline
$v = (3,3)$ & $\D^b(\PP^1_{2, 3, 3}) \times \D^b(\CC \vec{A}_{1})$ & $8 \cdot 1224720$ \\ \hline
$v = (3,4)$ & $\D^b(\PP^1_{2, 3, 4})$ & $46448640$ \\ \hline
\end{tabular}

\vspace{5pt}
\caption{Case of $A = (2, 3, 5)$}
\label{table : (2, 3, 5)}
\end{table}

Since $\chi_A = 1 / 30$, Theorem \ref{thm : recursion formula} and Table \ref{table : (2, 3, 5)} imply
\begin{eqnarray*}
e(\D^b(\PP^1_A)) & = & 30 (9^7 + 2 \cdot 7^8 + 8 \cdot 8^6 + 168 \cdot 3 \cdot 6^4 + 70 \cdot 5^3 \cdot 5^3 + \\
& & + 56 \cdot 2 \cdot 4^5 \cdot 4^2 + 28 \cdot 2^9 \cdot 3^4 \cdot 3 + 8 \cdot 2 \cdot 3^{12} + 2 \cdot 3^5 \cdot 5^7) \\ 
& & + 2 \cdot 105 \cdot 3^5 \cdot 5^7 + 3 (8 \cdot 30 \cdot 2^2 \cdot 5^5 + 21000000) \\
& & + 5 (56 \cdot 12 \cdot 2^2 \cdot 3^3 \cdot 4^2 + 28 \cdot 38880 \cdot 3 + 8 \cdot 1224720 + 46448640) \\
& = & 2551500000 \\
& = & \frac{9!}{2! 3! 5! \frac{1}{30}} 2^2 3^3 5^5.
\end{eqnarray*}
Hence, we obtain the statement for $A = (2, 3, 5)$.

We have finished the proof.
\qed
\section{Lyashko--Looijenga map of a Frobenius manifold}\label{sec : generalized LL map}\label{LL-Frob}

For a Frobenius manifold, one can define the Lyashko--Looijenga map.
The degree of the Lyashko--Looijenga map plays an important role in the paper.
In this section, we shall consider the Frobenius manifold associated with an orbifold projective line.
We also review the mirror symmetry among an affine cusp polynomial, an orbifold projective line and a generalized root system of affine ADE type.
We refer to \cite{T1} for more details of the mirror symmetry among them.

\subsection{Frobenius manifold for an affine cusp polynomial}\label{subsec : Frobenius manifold for an affine cusp polynomial}

In order to see some properties of the Frobenius manifold of an orbifold projective line via mirror symmetry, we first consider the Frobenius manifold associated with an affine cusp polynomial.
\begin{defn}
Let $A = (a_1, a_2, a_3)$ be a tuple satisfying $\chi_A > 0$.
A polynomial $f_A({\bf x}) \in \CC[x_1,x_2,x_3]$ given as
\begin{equation*}
f_A ({\bf x}) \coloneqq x_1^{a_1} + x_2^{a_2} + x_3^{a_3} - q^{-1} \cdot x_1 x_2 x_3
\end{equation*}
for a nonzero complex number $q \in \CC^\ast$ is called the {\em affine cusp polynomial} of type $A$.
\end{defn}
For an affine cusp polynomial $f_A ({\bf x})$, one can see that 
\[
\dim_\CC \CC[x_1, x_2, x_3] \left/ \Big( \frac{\p f_A}{\p x_1}, \frac{\p f_A}{\p x_2}, \frac{\p f_A}{\p x_3} \Big) \right. = \mu_A
\]
We shall consider the universal unfolding of $f_A$.
Define a complex manifold $M$ by 
\[
M \coloneqq \CC^{\mu_A - 1} \times \CC^\ast
\]
and denote by $({\bf s}, s_{\mu_A})$ the coordinates of $M$.
Let $p \colon \CC^3 \times M \longrightarrow M$ denote the natural projection map.

\begin{defn}\label{deform-affine cusp}
Define a holomorphic function $F_A \colon \CC^3 \times M \longrightarrow \CC$ by 
\begin{equation}
F_A({\bf x}; {\bf s}, s_{\mu_A}) \coloneqq x_1^{a_1} + x_2^{a_2} + x_3^{a_3} - s_{\mu_A}^{-1} \cdot x_1 x_2 x_3 + s_1 \cdot 1 + \sum_{i = 1}^3 \sum_{j = 1}^{a_i - 1} s_{i, j} \cdot x_i^j.
\end{equation}
\end{defn}

Set $p_* \O_\C \coloneqq \O_{M} [x_1, x_2, x_3] \left/ \Big( \frac{\p F_A}{\p x_1}, \frac{\p F_A}{\p x_2}, \frac{\p F_A}{\p x_3} \Big) \right. $.
The sheaf $p_* \O_\C$ can be thought of as the direct image of the sheaf of relative {\em algebraic} functions on the relative critical set $\C$ of $F_A$ with respect to the projection $p \colon \CC^3 \times M \longrightarrow M$.
Then, the pair $(F_A, p)$ forms the {\em universal unfolding} of $f_A$.
Namely, the following proposition holds.
\begin{prop}[{\cite[Proposition 2.5]{IST}}]\label{prop : universal unfolding}
The function $F_A$ satisfies the following conditions:
\begin{enumerate}
\item $F_A({\bf x}; {\bf 0}, q) = f_A({\bf x})$.
\item The $\O_M$-homomorphism called the Kodaira--Spencer map defined as
\begin{equation*}
\T_M \longrightarrow p_* \O_\C, \quad \delta \mapsto \delta F_A,
\end{equation*}
is an isomorphism.
\qed
\end{enumerate}
\end{prop}

We shall denote by $\circ$ the induced product structure on $\T_M$ by the Kodaira--Spencer map.
Namely, for $\delta, \delta' \in \T_M$, we have$(\delta \circ \delta') F_A = \delta F_A \cdot \delta' F_A$ in $p_* \O_\C$.
Define vector fields $e, E_A \in \Gamma(M, \T_M)$ by 
\[
e \coloneqq \dfrac{\p}{\p s_1}, \quad E_A \coloneqq s_1 \frac{\p}{\p s_1} + \sum_{i = 1}^3 \sum_{j = 1}^{a_i - 1} \frac{a_i - j}{a_i} s_{(i, j)} \frac{\p}{\p s_{(i, j)}} + \chi_A s_{\mu_A} \frac{\p}{\p s_{\mu_A}}.
\]
Then, the vector field $e$ is the unit of the $\O_M$-product $\circ$ on $\T_M$.
The vector field $E_A$ is called the {\em Euler vector field}.

In order to define the notion of a primitive form, it is necessary to define a Saito structure associated with the universal unfolding $F_A$. 
This structure is given as a tuple consisting of the filtered de Rham cohomology group $\H_{F_A}$ (whose increasing filtration is denoted by $\H_{F_A}^{(k)}$ $(k\in\ZZ)$), the Gau\ss--Manin connection $\nabla$ on $\H_{F_A}$ and the higher residue pairing $K_{F_A}$ on $\H_{F_A}$. 
In the paper, we omit the details about those objects and refer the interested reader to \cite{SaTa}. 
\begin{prop}[{\cite[Theorem 3.1]{IST}}]\label{prop : primitive form}
The element $\zeta_A \coloneqq [s_{\mu_A}^{-1} d x_1 \wedge d x_2 \wedge d x_3] \in \H_{F_A}^{(0)}$ is a primitive form for the tuple $(\H_{F_A}^{(0)}, \nabla, K_{F_A})$ with the minimal exponent $r = 1$.
\qed
\end{prop}

The higher residue pairing $K_{F_A}$ and the primitive form $\zeta_A$ induce a symmetric non-degenerate $\O_M$-bilinear form $\eta_{\zeta_A}:\T_M \times \T_M \longrightarrow \O_M$ defined by 
\begin{equation*}
\eta_{\zeta_A} (\delta, \delta') \coloneqq K_{F_A} \left( r^{(0)}(u \nabla_{\delta} \zeta_A),r^{(0)}(u \nabla_{\delta'} \zeta_A) \right), \quad \delta, \delta' \in \T_M ,
\end{equation*}
where $r^{(0)}$ is the natural map truncating $\H_{F_A}^{(-1)}$ component inside $\H_{F_A}^{(0)}$.

Let $M_{(f_A, \zeta_A)}$ denote the tuple $(M, \eta_{\zeta_A}, \circ, e, E_A)$.
\begin{prop}[{cf.~\cite{SaTa, D}}]\label{prop : Frobenius manifold}
The tuple $M_{(f_A, \zeta_A)} = (M, \eta_{\zeta_A}, \circ, e, E_A)$ is a Frobenius manifold of rank $\mu_A$ and dimension $1$. 
Namely, it satisfies the following conditions:
\begin{enumerate}
\item The product $\circ$ is self-adjoint with respect to $\eta_{\zeta_A}$: 
\begin{equation*}
\eta_{\zeta_A} (\delta \circ \delta', \delta'') = \eta_{\zeta_A} (\delta, \delta' \circ \delta''), \quad
\delta, \delta', \delta'' \in \T_M. 
\end{equation*} 
\item The Levi--Civita connection $\ns \colon \T_M \otimes_{\O_M} \T_M \longrightarrow \T_M$ with respect to $\eta_{\zeta_A}$ is
flat: 
\begin{equation*}
[\ns_\delta, \ns_{\delta'}] = \ns_{[\delta, \delta']}, \quad \delta, \delta' \in \T_M.
\end{equation*}
\item The $\O_M$-linear morphism $C_\delta \colon \T_M \longrightarrow \T_M, ~ \delta \in \T_M$ defined by $C_\delta(-) \coloneqq \delta \circ -$ satisfies 
\begin{equation*}
\ns_{\delta} (C_{\delta'} \delta'') - C_{\delta'} (\ns_{\delta} \delta'') - C_{{\nabla} \hspace{-1.0mm} \raisebox{0.3mm}{\text{\fontsize{6pt}{0cm}\selectfont{\bf /}}}_{\delta} \delta'} \delta'' = 
\ns_{\delta'} (C_{\delta} \delta'') - C_{\delta} (\ns_{\delta'} \delta'') - C_{{\nabla}\hspace{-1.0mm} \raisebox{0.3mm}{\text{\fontsize{6pt}{0cm}\selectfont{\bf /}}}_{\delta'} \delta} \delta''.
\end{equation*}
\item The unit vector field $e$ is $\ns$-flat: 
\begin{equation*}
\ns e = 0.
\end{equation*} 
\item The bilinear form $\eta_{\zeta_A}$ and the product $\circ$ are homogeneous of degree $1$ with respect to the Lie derivative of the Euler vector field $E_A$: 
\begin{equation*}
{\rm Lie}_{E_A} (\eta_{\zeta_A}) = \eta_{\zeta_A}, \quad 
{\rm Lie}_{E_A} (\circ) = \circ.
\end{equation*}
\end{enumerate}
\qed
\end{prop}

On the other hand, one can construct another Frobenius manifold from the viewpoint of generalized root systems. 
It is known by \cite{STW} that, for a Dynkin quiver or an extended Dynkin quiver, a generalized root system is obtained by the derived category. 
Based on the works \cite{S-K, SYS}, Dubrovin \cite{D} constructed a Frobenius manifold for the generalized root system of a Dynkin quiver by the Weyl group invariant theory.
As an analogue of the work \cite{D}, Dubrovin--Zhang generalized the construction to the cases of extended Dynkin quivers.
By the invariant theory of the extended Weyl group $\widehat{W}_A$, the Frobenius structure is given on the quotient space $\widehat{\h}_A / \widehat{W}_A$, where $\widehat{\h}_A = \h_A \times \CC$ is the extension of the $(\mu_A - 1)$-dimensional complex vector space $\h_A$.
The complex manifold $\widehat{\h}_A / \widehat{W}_A$ is isomorphic to $M$.
We denote by $M_{\widehat{W}_A}$ the Frobenius manifold for an extended Dynkin quiver $Q_A$.
Then, the next proposition is known.
\begin{prop}[{\cite{DZ}, \cite[Corollary 7.3]{ShTa}}]\label{prop : isomorphism of Frobenius manifolds 1}
The Frobenius manifold $M_{\widehat{W}_A}$ is isomorphic to $M_{(f_A, \zeta_A)}$.
\qed
\end{prop}

\subsection{Frobenius manifold for an orbifold projective line}

Let $\PP_A^1$ be an orbifold projective line with $\chi_A>0$.
We shall consider a Frobenius structure on $M = \CC^{\mu_A - 1} \times \CC^\ast$ constructed from the orbifold Gromov--Witten theory for $\PP_A^1$.
One can define the orbifold cohomology $H_{\rm orb}^* (\PP_A^1)$ of $\PP_A^1$ equipped with the orbifold Poincar\'{e} pairing 
\[
(-,-)_{\rm orb} \colon H_{\rm orb}^* (\PP_A^1) \times H_{\rm orb}^* (\PP_A^1) \longrightarrow \CC, 
\]
which is a non-degenerate symmetric $\CC$-bilinear form.
The orbifold cohomology $H_{\rm orb}^* (\PP_A^1)$ is a $\mu_A$-dimensional $\CC$-vector space such that there is a $\CC$-basis 
\begin{equation}\label{eq : basis of the orbifold cohomology}
\{ \phi_1, \phi_{(1,1)}, \dots, \phi_{(1, a_1 - 1)}, \phi_{(2, 1)}, \dots, \phi_{(2, a_2 - 1)}, \phi_{(3, 1)}, \dots, \phi_{(3, a_3 - 1)}, \phi_{\mu_A} \}
\end{equation}
satisfying 
\begin{equation*}
H_{\rm orb}^0 (\PP_A^1) \cong \CC \phi_1, \quad H_{\rm orb}^2 (\PP_A^1) \cong \CC \phi_{\mu_A},
\end{equation*}
\begin{equation*}
\phi_{(i, j)} \in H_{\rm orb}^{\frac{2 j}{a_i}} (\PP_A^1), \quad i = 1, 2, 3, ~ j = 1, \dots, a_i - 1.
\end{equation*}
and 
\begin{equation*}
(\phi_1, \phi_{\mu_A})_{\rm orb} = 1, \quad 
(\phi_1, \phi_{(i, j)})_{\rm orb} = 0, \quad 
(\phi_{(i, j)}, \phi_{\mu_A})_{\rm orb} = 0, 
\end{equation*}
\begin{equation*} 
(\phi_{(i, j)}, \phi_{(i', j')})_{\rm orb} = 
\begin{cases}
\dfrac{1}{a_i}, & (i = i', ~ j = a_i - j'), \\
0, & (\text{otherwise}).
\end{cases}
\end{equation*}

Let $(t_1, \dots, t_{\mu_A})$ denote the dual coordinates of the $\CC$-basis \eqref{eq : basis of the orbifold cohomology}.
Then, 
\[
(t_1, t_{(1,1)}, \dots, t_{(1, a_1 - 1)}, t_{(2, 1)}, \dots, t_{(2, a_2 - 1)}, t_{(3, 1)}, \dots, t_{(3, a_3 - 1)}, e^{t_{\mu_A}})
\]
forms a coordinate system on $M$.
Note that there is an $\O_{M}$-isomorphism between the holomorphic tangent sheaf $\T_{M}$ and the sheaf $H^*_{\rm orb} (\PP_A^1) \otimes_\CC \O_{M}$:
\begin{equation}\label{eq : tangent sheaf}
H^*_{\rm orb} (\PP_A^1) \otimes_\CC \O_{M} \cong \T_{M}, \quad \phi_k \mapsto \frac{\p}{\p t_k}.
\end{equation}
Define holomorphic vector fields $e$ and $E_A$ on $\T_M$ by 
\[
e \coloneqq \dfrac{\p}{\p t_1}, \quad E_A \coloneqq t_1 \frac{\p}{\p t_1} + \sum_{i = 1}^3 \sum_{j = 1}^{a_i - 1} \frac{a_i - j}{a_i} t_{(i, j)} \frac{\p}{\p t_{(i, j)}} + \chi_A \frac{\p}{\p t_{\mu_A}}.
\]
Denote by $\eta \colon \T_M \times \T_M \longrightarrow \O_M$ the non-degenerate symmetric $\O_M$-bilinear form induced by the orbifold Poincar\'{e} pairing $(-,-)_{\rm orb}$ and the $\O_M$-isomorphism \eqref{eq : tangent sheaf}.

For $n \in \ZZ_{\ge 0}$ and $\beta \in H_2(\PP_A^1, \ZZ)$, the moduli stack $\overline{\M}_{0, n}(\PP_A^1, \beta)$ of orbifold (twisted) stable maps of genus $0$ with $n$-marked points of degree $\beta$ is defined by Chen--Ruan \cite{CR} and Abramovich--Graber--Vistoli \cite{AGV}.
It is shown that there exists a virtual fundamental class $[\overline{\M}_{0, n}(\PP_A^1, \beta)]^{\rm vir}$ and Gromov--Witten invariants of genus $0$ with $n$-marked points of degree $\beta$ are defined as
\begin{equation*}
\langle \alpha_1, \dots, \alpha_n \rangle_{0, n, \beta}^{\PP_A^1} \coloneqq
\int_{[\overline{\M}_{0, n}(\PP_A^1, \beta)]^{\rm vir}} ev_1^\ast \alpha_1 \wedge \dots \wedge ev_n^\ast \alpha_n
\end{equation*}
for $\alpha_1, \dots, \alpha_n \in H_{\rm orb}^* (\PP_A^1, \QQ)$, where $ev_i^\ast \colon H_{\rm orb}^* (\PP_A^1, \QQ) \longrightarrow H^* (\overline{\M}_{0, n}(\PP_A^1, \beta), \QQ)$ denotes the induced homomorphism by the evaluation map at $i$-th marked point.
The generating function 
\begin{equation*}
\F_0^{\PP_A^1} \coloneqq 
\sum_{\substack{n \in \ZZ_{\ge 0} \\ \beta \in H_2(\PP_A^1, \ZZ)}} \frac{1}{n !} \langle \mathbf{t}, \dots, \mathbf{t} \rangle_{0, n, \beta}^{\PP_A^1}, 
\quad \mathbf{t} = \sum_{i = 1}^{\mu_A} t_i \phi_i
\end{equation*}
is called the {\em genus zero Gromov--Witten potential}.
The genus zero Gromov--Witten potential defines the {\em quantum product} on $H_{\rm orb}^* (\PP_A^1) \otimes_\CC \O_M$, which yields a product structure $\circ \colon \T_M \times \T_M \longrightarrow \T_M$ by the $\O_M$-isomorphism \eqref{eq : tangent sheaf}.
Due to Chen--Ruan \cite{CR} and Abramovich--Graber--Vistoli \cite{AGV}, the tuple $M_{\PP_A^1} \coloneqq (M, \eta, \circ, e, E_A)$ is a (formal) Frobenius manifold.

Milanov--Tseng and Rossi proved the mirror symmetry for an orbifold projective line with $\chi_A>0$.
Combined with Proposition \ref{prop : isomorphism of Frobenius manifolds 1}, the following holds.
\begin{prop}[Classical mirror symmetry \cite{MT, Ro} cf.~\cite{IST}]\label{cor : CMS}
There exists an isomorphism of Frobenius manifolds 
\begin{equation*}
M_{\PP_A^1} \cong M_{(f_A, \zeta_A)} \cong M_{\widehat{W}_A}.
\end{equation*}
\qed
\end{prop}
\begin{rem}[Homological mirror symmetry]
For an affine cusp polynomial $f_A$, one can consider a triangulated category $\D^b\Fuk^\to(f_A)$, which is called the derived {\em Fukaya--Seidel} category \cite{S-P}.
It is known by \cite{GL1, S-P, S-D, AKO, T2} that there exist equivalences 
\[
\D^b\coh(\PP_A^1) \cong \D^b\Fuk^\to(f_A) \cong \D^b\mod(\CC Q_A).
\]
These equivalences are called the {\em homological mirror symmetry}.
See \cite{T1} for more details about the mirror symmetry.
\end{rem}

It was proved by \cite{K} that the Frobenius manifold $M_{\PP_A^1}$ is semisimple.
We can also obtain the fact as a consequence of Proposition \ref{cor : CMS}.
\begin{prop}[{\cite{K}}]
The Frobenius manifold $M_{\PP_A^1}$ is semi-simple.
In particular, for a general point $t \in M_{\PP_A^1}$, there exists a local coordinate system $(u_1, \dots, u_{\mu_A})$ such that 
\begin{equation*}
e = \frac{\p}{\p u_1} + \dots + \frac{\p}{\p u_{\mu_A}}, \quad 
E_A = u_1 \frac{\p}{\p u_1} + \dots + u_{\mu_A} \frac{\p}{\p u_{\mu_A}}, 
\end{equation*}
\begin{equation*}
\frac{\p}{\p u_i} \circ \frac{\p}{\p u_j} = \delta_{ij} \frac{\p}{\p u_i}, \quad i, j = 1, \dots, \mu_A, 
\end{equation*}
where $\delta_{ij}$ is the Kronecker's delta.
\qed
\end{prop}
The local coordinate system $(u_1, \dots, u_{\mu_A})$ is called the {\em canonical coordinate system}.
The canonical coordinate system is uniquely determined up to a permutation of indices.

\begin{defn}[{\cite[Section 3]{DZ}}]
For a semi-simple Frobenius manifold $M$ of rank $\mu$ with the Euler vector field $E$, more generally, for a massive $F$-manifold (see {\cite[Section 3.5]{He}}), 
the {\em Lyashko--Looijenga map} $LL \colon M \longrightarrow \CC^\mu$ is defined as
\begin{equation*}
LL (t) \coloneqq (b_1, \dots, b_{\mu}),
\end{equation*}
where $b_1, \dots, b_{\mu_A}$ are coefficients of the characteristic polynomial of $C_E$:
\begin{equation*}
\det (C_{E} - w) = \prod_{i = 1}^{\mu} (u_i - w) = (-1)^\mu (w^{\mu}  + b_1 w^{\mu- 1} + \dots + b_{\mu}).
\end{equation*}
\end{defn}
It is known that the Lyashko--Looijenga map is locally bi-holomorphic on the dense open subset of $M$, complement of the bifurcation set (see \cite[Theorem~3.19]{He} for the precise statement). 

In the case of $M_{\PP_A^1} \cong M_{(f_A, \zeta_A)} \cong M_{\widehat{W}_A}$,  
the Lyashko--Looijenga map $LL$ is a ramified covering map and its {\em degree} $\deg LL$ is given explicitly as follows.
Since the degrees of parameters in \eqref{deform-affine cusp} 
$\displaystyle \deg s_{1}=\deg t_{1}=1$, $\displaystyle \deg s_{(i,j)}=\deg t_{(i,j)}=\frac{a_{i}-j}{a_{i}}$
and $\displaystyle \deg s_{\mu_{A}}=\deg e^{t_{s_{\mu_{A}}}}=\chi_{A}$ are all positive, 
the Frobenius potential for $M_{\PP_A^1} \cong M_{(f_A, \zeta_A)} \cong M_{\widehat{W}_A}$ becomes 
a weighted  homogeneous polynomial with respect to the variables $t_{1},t_{(1,1)},\dots, t_{(3,a_{3}-1)}, e^{t_{\mu_{A}}}$.
This implies all entries $b_{i}$ of the Lyashko--Looijenga map are also 
weighted homogeneous polynomials with respect to the variables $t_{1},t_{(1,1)},\dots, t_{(3,a_{3}-1)}, e^{t_{\mu_{A}}}$ as in the argument 
of \cite[Section 3]{DZ}.
In particular, the weighted degree of $b_{i}$ is just $i$ since the degrees of canonical coordinates $u_{i}$ are just $1$.
To summarize, we obtain the following  
\begin{prop}[{\cite[Section 3]{DZ}}]\label{prop : deg LL}
We have
\begin{equation}
\deg LL = \dfrac{\mu_A !}{\displaystyle \chi_A \prod_{i = 1}^3 \prod_{j = 1}^{a_i - 1} \dfrac{a_i - j}{a_i}}=\frac{\mu_A!}{a_1!a_2!a_3!\chi_A}a_1^{a_1}a_2^{a_2}a_3^{a_3} .
\end{equation}
\qed
\end{prop}

As a direct consequence of Theorem \ref{thm : main}, the number $e(\D^b(\PP_A^1))$ is equal to the degree of the Lyashko--Looijenga map for the Frobenius manifold $M_{\PP_A^1}$.
\begin{cor}\label{cor: equality to deg LL}
We have 
\[
e(\D^b(\PP_A^1)) = \deg LL.
\]
\qed
\end{cor}

According to the third named author, Corollary~\ref{cor: equality to deg LL} extends to the case when $\chi_A=0$ 
after appropriately adapting the definition to the situation where $\deg LL$ is the one calculated by \cite{HR} for simple elliptic singularities in the Legendre normal forms with $A$
as the Gabrielov numbers. The detail will be reported elsewhere \cite{TZ}.

\end{document}